\documentclass[10pt]{article}

\usepackage{a4wide}
\usepackage{amssymb}
\usepackage{amsfonts}
\usepackage{amsmath}
\input xy
\xyoption{arrow} \xyoption{matrix}

\date{}

\newtheorem{proposition}{Proposition}[section]
\newtheorem{theorem}[proposition]{Theorem}
\newtheorem{lemma}[proposition]{Lemma}

\newtheorem{corollary}[proposition]{Corollary}

\def\Hom{{\rm Hom}}
\def\der{\partial }

\def\nFM0{{\nu }_{F,M_0}}
\def\nFN0{{\nu }_{F,N_0}}
\def\nGN0{{\nu }_{G,N_0}}

\def\N0{ {\bf N}_0 }

\def\g{\gamma}
\def\v{\varphi}
\def\ra{\rightarrow}

\def\Xpm{X^{\pm }}

\def\s{\sigma}
\def\Z{\mathbb{Z}}

\def\l1{{\lambda}_1}

\def\a{\alpha}
\def\a0{ {\alpha }_0}
\def\a1{ {\alpha }_1}

\def\l{\lambda}
\def\o{\omega}

\def\nFGM0{{\nu }_{F,G,M_0}}

\def\nFN0{{\nu}_{F,N_0}}

\def\sm{{\sigma}^m}

\def\sm1{{\sigma}^{-1}}

\def\smtp1{{\sigma}^{-t+1}}

\def\o{\omega }
\def\S1{S^{-1}}

\def\Xpm1{X^{\pm 1}_1}

\def\sPM1{{\sigma }^{\pm 1}}
\def\sMP1{{\sigma }^{\mp 1 }}

\def\d{\delta}

\def\di{{\rm d.ind}}

\def\L{\Lambda}
\def\O{\Omega}

\def\CA{{\cal A}}

\def\CD{{\cal D}}

\def\Ytm1{Y^{t-1}}
\def\Yim1{Y^{i-1}}

\def\CM{{\cal M}}

\def\CF{{\cal F}}
\def\CG{{\cal G}}

\def\Aut{{\rm Aut}}

\def\Der{{\rm Der }}
\def\ad{{\rm ad }}
\def\dim{{\rm dim }}

\def\ker{ {\rm ker } }

\def\CJ{ {\cal J}}

\def\D{ \Delta }

\def\SL2Z{ {\rm SL}_2({\bf Z}) }

\def\Gp1{ G^{1 , 1 } }
\def\P11{ P^{-1 , 1 } }
\def\Pp1{ P^{1 , 1 } }

\def\nCLsr{{}^\nu\kern-2pt {\cal L}^{\sigma , \rho  }}
\def\nP{{}^\nu \kern-2pt P}
\def\nL{{}^\nu\kern-2pt L}
\def\nLL{{}^\nu\kern-2pt \Lambda}
\def\nPsr{{}^\nu\kern-2pt P^{\sigma , \rho  }}
\def\nLsr{{}^\nu\kern-2pt L^{\sigma , \rho  }}
\def\nuCL{{}^\nu\kern-2pt  {\cal L}}
\def\nCLsr{{}^\nu\kern-2pt {\cal L}^{\sigma , \rho  }}
\def\nCL1m{{}^\nu\kern-2pt {\cal L}^{-1 , 1  }}
\def\x1nu{x^\frac{1}{\nu}}
\def\xm1nu{x^{-\frac{1}{\nu}}}

\def\ra{\rightarrow }

\def\CB{{\cal B}}

\def\CI{{\cal I}}

\def\CT{{\cal T}}

\def\nAM0{{\nu }_{{\cal A},M_0}}
\def\nAN0{{\nu }_{{\cal A},N_0}}

\def\End{ {\rm End }}
\def\Der{ {\rm Der }}
\def\CJ{ {\cal J }}

\def\ad{ {\rm ad }}

\def\ga{\mathfrak{a}}
\def\gb{\mathfrak{b}}

\def\gm{\mathfrak{m}}

\def\SL{{\rm SL}}

\def\Hom{{\rm Hom}}

\def\di!{\frac{\der^i}{i!}}
\def\dik!{\frac{\der^k_i}{k!}}

\def\ord{{\rm ord}}

\def\id{{\rm id}}

\def\N{\mathbb{N}}

\def\0{\overline{0}}
\def\1{\overline{1}}

\def\Ln1{\L_{n,\overline{1}}}

\def\a1{a_{\overline{1}}}

\def\S{\Sigma}

\def\CU{{\cal U}}

\def\vn1{\overrightarrow{n-1}}

\def\Sh{{\rm Sh}}

\def\im{{\rm im}}

\def\Inn{{\rm Inn}}

\def\mJ{\mathbb{J}}
\def\mI{\mathbb{I}}

\def\ann{{\rm ann}}

\def\Cen{{\rm Cen}}

\def\mF{\mathbb{F}}

\def\mT{\mathbb{T}}

\def\mE{\mathbb{E}}

\def\K1{{\rm K}_1}

\def\hmI1{\widehat{\mI_1}}
\def\tmI1{\widetilde{\mI_1}}
\def\tmJ1{\widetilde{\mJ_1}}
\def\hB1{\widehat{B_1}}
\def\hCB1{\widehat{\CB_1}}

\def\mW{\mathrm{W}}

\def\ggu{\mathfrak{u}}
\def\Xai{X_{\alpha , i}}
\def\Xbj{X_{\beta , j}}
\def\udim{{\rm u.dim}}

\def\Fix{{\rm Fix}}
\def\TAut{{\rm TAut}}
\def\UAut{{\rm UAut}}
\def\sh{{\rm sh}}

\def\mJ{\mathbb{J}}
\def\res{{\rm res}}

\def\LN{{\rm LN}}

\begin{document}

\author{V. V. \  Bavula   
}

\title{The groups of automorphisms of the Lie algebras of triangular polynomial derivations }

\maketitle

\begin{abstract}

The group of automorphisms $G_n$ of the Lie algebra $\ggu_n$ of
triangular polynomial derivations of the polynomial algebra $P_n =
K[x_1,\ldots , x_n]$  is found ($n\geq 2$), it is isomorphic to an
iterated semi-direct  product
 $$\mT^n\ltimes (\UAut_K(P_n)_n\rtimes( \mF_n' \times  \mE_n ))  $$
 where $\mT^n$ is an algebraic  $n$-dimensional torus,   $\UAut_K(P_n)_n$ is an explicit factor group of the group $\UAut_K(P_n)$ of triangular polynomial automorphisms, $\mF_n'$ and $\mE_n$ are explicit groups that are isomorphic respectively to the groups $\mI$ and $\mJ^{n-2}$ where
    $\mI := (1+t^2K[[t]], \cdot )\simeq K^{\N}$ and
  $\mJ := (tK[[t]], +)\simeq K^\N$. It is shown that the adjoint group of automorphisms of the Lie algebra $\ggu_n $ is equal to the group $\UAut_K(P_n)_n$.

$\noindent $

{\em Key Words: Group of automorphisms, Lie algebra, triangular
polynomial derivations, automorphism,  locally nilpotent
derivation. }

 {\em Mathematics subject classification
2010:  17B40, 17B66,  17B65, 17B30.}

$${\bf Contents}$$
\begin{enumerate}
\item Introduction.
 \item The structure of the group of automorphisms of the  Lie algebra $\ggu_n$.
\item  The group of automorphisms of the  Lie algebra $\ggu_n$.
\item The group of automorphism of the   Lie algebra $\ggu_n$ is an iterated semi-direct  product.
\item The canonical decomposition for an  automorphism of the   Lie algebra $\ggu_n$.
\item The adjoint group of automorphisms  of the Lie algebra $\ggu_n$.

\end{enumerate}
\end{abstract}


\section{Introduction}

Throughout, 
 module means
a left module;
 $\N :=\{0, 1, \ldots \}$ is the set of natural numbers; $K$ is a
field of characteristic zero and  $K^*$ is its group of units;
$P_n:= K[x_1, \ldots , x_n]=\bigoplus_{\alpha \in \N^n}
Kx^{\alpha}$ is a polynomial algebra over $K$ where
$x^{\alpha}:=x_1^{\alpha_1}\cdots x_n^{\alpha_n}$;
$\der_1:=\frac{\der}{\der x_1}, \ldots , \der_n:=\frac{\der}{\der
x_n}$ are the partial derivatives ($K$-linear derivations) of
$P_n$;  $\Aut_K(P_n)$ is the group of automorphisms of the polynomial algebra $P_n$; $\Der_K(P_n) =\bigoplus_{i=1}^nP_n\der_i$ is the Lie
algebra of $K$-derivations of $P_n$; $A_n:= K \langle x_1, \ldots
, x_n , \der_1, \ldots , \der_n\rangle  =\bigoplus_{\alpha , \beta
\in \N^n} Kx^\alpha \der^\beta$  is  the $n$'th {\em Weyl
algebra}; for each natural number $n\geq 2$,
$$\ggu_n :=
K\der_1+P_1\der_2+\cdots +P_{n-1}\der_n$$ is the  {\em Lie algebra
of triangular polynomial derivations} (it is a Lie subalgebra of
the Lie algebra $\Der_K(P_n)$) and $G_n:= \Aut_K(\ggu_n)$ is its
group of automorphisms;
$\d_1:=\ad (\der_1), \ldots , \d_n:=\ad (\der_n)$ are the inner derivations of the Lie algebra $\ggu_n$ determined by the elements $\der_1, \ldots , \der_n$ (where $\ad (a)(b):=[a,b]$).

{\bf The group of automorphisms $G_n$ of the Lie algebra $\ggu_n$}.
The aim of the paper is to find the group $G_n$ (Theorem \ref{28Feb12}) and its explicit generators.
\begin{itemize}
\item {\rm (Theorem \ref{28Feb12})}
{\em Let $\mI := (1+t^2K[[t]], \cdot)$ and $\mJ:= (tK[[t]], +)$. Then for all} $n\geq 2$,
\begin{enumerate}
\item $G_n = \mT^n \ltimes (\UAut_K(P_n)_n \rtimes ( \mF_n'\times \mE_n))$.
\item $ G_n\simeq \mT^n \ltimes (\UAut_K(P_n)_n \rtimes (\mI \times \mJ^{n-2})$.
    \end{enumerate}
\end{itemize}
The group $\mT^n$ is an algebraic  $n$-dimensional torus,
 $\UAut_K(P_n)_n:=\UAut_K(P_n)/\sh_n$ is the factor group of the group of triangular polynomial automorphisms
$$\UAut_K(P_n):=\{ \s \in \Aut_K(P_n)\, | \, \s (x_i) = x_i+a_i, \;  a_i\in P_{i-1}\; {\rm for}\; i=1, \ldots , n\}$$
modulo its normal subgroup
$$ \sh_n := \{ \s \in \Aut_K(P_n)\, | \, \s (x_i) = x_i, \; i=1, \ldots , n-1; \; \s (x_n) = x_n +\l , \; \l \in K\}, $$
 $\mF_n'\simeq \mI$ and $\mE_n\simeq \mJ^{n-2}$ are explicit subgroups of $G_n$ (see below and Section \ref{GAUTUN}). The group $G_n$ is made up  of two parts: the `obvious' one, $\mT^n \ltimes \UAut_K(P_n)_n$,  and the `non-obvious' one -- $\mF_n'\times \mE_n\simeq \mI \times \mJ^{n-2}$ -- which is a much  more massive group than the group $\mT^n \ltimes \UAut_K(P_n)_n$.

{\bf The key ideas and the strategy of finding the group $G_n$}.  A group $G=G_1\times_{ex} G_2$ is an {\em exact} product of its two subgroups $G_1$ and $G_2$ if every element $g$ of the group $G$ is a unique product $g_1g_2$ for some (unique)  elements $g_1\in G_1$ and $g_2\in G_2$.
 The strategy of finding the group $G_n$ is a (rather long) `refining process' which is done in  Sections \ref{GALAUN}--\ref{GnITER}. It consists of several steps. On each step the group $G_n$ is presented as an exact or semi-direct product of several of its subgroups.  Some of  these subgroups are explicit groups  and the other are defined in abstract terms (i.e., they satisfy certain properties, elements of which satisfy certain equations). Every successive step is a `refinement' of its predecessor in the sense that `abstract' subgroups are presented  as explicit sets of automorphisms (i.e., the solutions are found to the defining equations of the subgroups).

In Section \ref{GALAUN}, the first step is done on the way  of
finding  the group $G_n$. In Section \ref{GALAUN}, several
important subgroups of the group $G_n$ are introduced.  These
include the group $\TAut_K(P_n)_n$  and its subgroup $\CT_n$ of
triangular polynomial automorphisms with all constant terms being
equal to zero, $$\CT_n:=\{ \s \in \Aut_K(P_n)\, | \, \s (x_1)=x_1,
\s(x_i) = x_i+a_i\; {\rm where}\;\; a_i\in (x_1, \ldots ,
x_{i-1}), i=2, \ldots , n\}$$ where $ (x_1, \ldots , x_{i-1})$ is
the maximal ideal of the polynomial algebra $P_{i-1}$ generated by
the elements $x_1, \ldots , x_{i-1}$;  and the group
$$\Sh_{n-1}:= \{ \s \in \Aut_K(P_n)\, | \, \s (x_1)=x_1+\l_1, \ldots , \s (x_{n-1})= x_{n-1}+\l_{n-1}, \s (x_n) = x_n\;{\rm  where}\;\;  \l_i\in K\}.$$ The most important subgroup of the group $G_n$ is
$$ \CF_n := \{ \s \in G_n \, | \, \s (\der_1) = \der_1, \ldots , \s (\der_n) = \der_n \},$$
as Theorem \ref{26Jan12} demonstrates.

\begin{itemize}
\item {\rm (Theorem \ref{26Jan12})}
\begin{enumerate}
\item $G_n=\TAut_K(P_n)_n \CF_n = \CF_n \TAut_K(P_n)_n$ {\em and}  $\TAut_K(P_n)_n \cap \CF_n = \Sh_{n-1}$.
\item $G_n=\mT^n\ltimes (\CT_n\times_{ex} \CF_n )=\mT^n\ltimes ( \CF_n \times_{ex} \CT_n)$.
\end{enumerate}
\end{itemize}
As the groups $\mT^n$ and $\CT_n$ are explicit groups, the problem of finding the group $G_n$ boils down to the problem of finding the group $\CF_n$. This is done in Section \ref{GAUTUN}.
\begin{itemize}
\item {\rm (Theorem \ref{11Feb12})} $\;\;\CF_n=\Sh_{n-2}\times \mF_n \times \mE_n$,
\end{itemize}
where \begin{eqnarray*}
 \Sh_{n-2}&=&\{ \s \in \Aut_K(P_n)\, | \, \s (x_1) = x_1+\l_1, \ldots ,\s (x_{n-2}) = x_{n-2}+\l_{n-2},\\
 & & \s (x_{n-1}) = x_{n-1}, \s (x_n) = x_n,  \; {\rm where} \; \l_i\in K \};  \\
 \mF_n &=& \{ f\in 1+\der_{n-1}K[[\der_{n-1}]]\; | \; f(p_i\der_i):=\begin{cases}
p_i\der_i & \text{if }i=1, \ldots , n-1,\\
f(p_n)\der_n & \text{if }i=n,
\end{cases} \\
 & & {\rm where } \; p_i\in P_{i-1}, i=1, \ldots , n \};\\
\mE_n &=&
\begin{cases}
\{ e \} & \text{if }n=2,\\
\prod_{j=2}^{n-1}\mE_{n,j}& \text{if }n\geq 3,\\
\end{cases} \\
\mE_{n,j}&=& \{ e_j'\in \der_{j-1}K[[\der_{j-1}]]\; | \;  e_j'(p_i\der_i) :=
\begin{cases}
p_j\der_j+e_j'(p_j)\der_n & \text{if }i=j,\\
p_i\der_i & \text{if }i\neq j.\\
\end{cases}
\end{eqnarray*}
As  a corollary, the group $G_n$ is presented as an exact product of its explicit subgroups.

\begin{itemize}
\item {\rm (Theorem \ref{5Feb12})}
 {\em Let $\mI = (1+t^2K[[t]], \cdot)$ and $\mJ = (tK[[t]], +)$. Then for all} $n\geq 2$,
\begin{enumerate}
\item $G_n = \mT^n \ltimes (\CT_n \times_{ex} (\Sh_{n-2}\times \mF_n\times \mE_n))=\TAut_K(P_n)_n\times_{ex} (\mF_n'\times \mE_n)$,
\item $ G_n\simeq \TAut_K(P_n)_n\times_{ex} (\mI \times \mJ^{n-2})$.
\end{enumerate}
\end{itemize}

In Section \ref{GnITER}, the explicit form of the groups $\mT^n$, $\UAut_K(P_n)_n$, $\mF_n$ and $\mE_n$ allows us to establish commutation relations between elements of these groups (Lemma \ref{a28Feb12} and Lemma \ref{b28Feb12}). From which we deduce that the group $\UAut_K(P_n)_n$ is a {\em normal} subgroup of the group $G_n$. In combination with Theorem \ref{5Feb12}.(1), this fact yields the main result of the paper
 $G_n = \mT^n \ltimes (\UAut_K(P_n)_n\rtimes (\mF_n'\times \mE_n))$ (Theorem \ref{28Feb12}.(1)), where $\mF_n'=1+\der_{n-1}^2K[[\der_{n-1}]]\subset \mF_n$.

 At the end of Section \ref{GnITER}, characterizations of the groups $\mF_n$ , $\mF_n'$ and $\mE_n$ are given in invariant terms (Proposition \ref{b12Feb12}).


{\bf The canonical decomposition for an automorphism of the Lie algebra $\ggu_n$}.

By Theorem \ref{5Feb12}.(1), every automorphism $\s \in G_n=\mT^n \ltimes (\CT_n\times_{ex}(\Sh_{n-2}\times \mF_n \times \mE_n))$ is the unique product
$$\s = t \tau s fe'\;\; {\rm  where}\;\; t\in \mT^n,\;  \tau \in \CT_n,\; s\in \Sh_{n-2},\;  f\in \mF_n,\;  e'\in \mE_n .$$
 This product is called the {\em canonical decomposition}  of  the automorphism $\s\in G_n$. In Section \ref{CDEC}, for every automorphism $\s \in G_n$ explicit formulas are found (Theorem \ref{14Feb12}) for the automorphisms $t$, $\tau$, $s$, $f$ and $e'$ via the elements $\{ \s (s) \, | \, s\in S_n\}$ where the set $S_n: =\{ \der_1, x_1^j\der_2, \ldots , x_{i-1}^j\der_i, \ldots , x_{n-1}^j\der_n\, | \,$ $  j\in \N\}$ is a set of  generators for the Lie algebra $\ggu_n$.

{\bf The adjoint group $\CA (\ggu_n)$ of the Lie algebra $\ggu_n$}. For a Lie algebra $\CG$, the {\em adjoint group} $\CA (\CG )$ is the subgroup of the group of automorphisms $\Aut_K(\CG )$ of the Lie algebra $\CG$ generated by the automorphisms $e^\d :=\sum_{i\geq 0}\frac{\d^i}{i!}$ where $\d$ runs through the set of    locally nilpotent inner derivations of the Lie algebra $\CG$. All the inner derivations of  the Lie algebra $\ggu_n$ are locally nilpotent derivations \cite{Lie-Un-GEN}. In Section \ref{ADJGRP}, we prove that the adjoint group $\CA (\ggu_n)$ of the Lie algebra $\ggu_n$ is equal to the group $\UAut_K(P_n)_n$ (Theorem \ref{8Apr12}).


\section{The  Lie algebra
$\ggu_n$}\label{TLAUN}

In this section, for reader's convenience various results and  properties of the Lie algebras $\ggu_n$ are collected that are used in the rest of the paper. The details/proofs  can be found in \cite{Lie-Un-GEN}.
 Since  $ \ggu_n=\bigoplus_{i=1}^n\bigoplus_{\alpha \in
\N^{i-1}}Kx^\alpha \der_i$, the elements 
\begin{equation}\label{Xai}
\Xai := x^\alpha \der_i = x_1^{\alpha_1}\cdots
x_{i-1}^{\alpha_{i-1}}\der_i, \;\;\; i=1, \ldots , n;  \; \alpha =(\alpha_1, \ldots , \alpha_n)
\in \N^{i-1},
\end{equation}
form the $K$-basis $\CB_n$  for the Lie algebra $\ggu_n$. The basis
$\CB_n$  is called the {\em canonical basis} for $\ggu_n$.  For
all $1\leq i\leq j\leq n$, $\alpha \in \N^{i-1}$ and $\beta \in
\N^{j-1}$, 
\begin{equation}\label{Xai1}
[\Xai ,\Xbj]=\begin{cases}
0& \text{if }i=j,\\
\beta_iX_{\alpha +\beta -e_i, j}& \text{if }i<j,\\
\end{cases}
\end{equation}
where $e_1:=(1,0, \ldots , 0), \ldots , e_n:=(0,\ldots, 0, 1)$ is
the canonical free $\Z$-basis for the $\Z$-module $\Z^n$. The Lie
algebra $\ggu_n=\oplus_{i=1}^nP_{i-1}\der_i$ is the direct sum of
{\em abelian} (infinite dimensional when $i>1$)  Lie subalgebras $P_{i-1}\der_i$
(i.e.,  $[P_{i-1}\der_i, P_{i-1}\der_i]=0$) such that, for all $i<j$,
\begin{equation}\label{PdiPdj}
[P_{i-1}\der_i, P_{j-1}\der_j]=P_{j-1}\der_j.
\end{equation}
The Lie subalgebra $P_{i-1}\der_i$ has the structure of the left $P_{i-1}$-module and ${}_{P_{i-1}}(P_{i-1}\der_i)\simeq P_{i-1}$.  By (\ref{PdiPdj}), the Lie algebra $\ggu_n$
admits the finite strictly descending chain of ideals
\begin{equation}\label{seruni}
\ggu_{n,1}:=\ggu_n\supset \ggu_{n,2}\supset \cdots \supset
\ggu_{n,i}\supset \cdots \supset \ggu_{n,n}\supset \ggu_{n,n+1}:=0
\end{equation}
where $\ggu_{n,i}:=\sum_{j=i}^nP_{j-1}\der_j$ for $i=1, \ldots ,
n$. By (\ref{PdiPdj}), for all $i<j$, 
\begin{equation}\label{cuni}
[\ggu_{n,i}, \ggu_{n,j}]\subseteq
\begin{cases}
\ggu_{n,i+1}& \text{if }i=j,\\
\ggu_{n,j}& \text{if }i<j.
\end{cases}
\end{equation}
For all $i=1, \ldots , n$, there is the canonical isomorphism of
Lie algebras 
\begin{equation}\label{cuni1}
\ggu_i\simeq \ggu_n/\ggu_{n,i+1}, \;\; X_{\alpha , j}\mapsto
X_{\alpha , j}+\ggu_{n, i+1}.
\end{equation}
 In particular, $\ggu_{n-1}\simeq \ggu_n/P_{n-1}\der_n$. The polynomial
algebra $P_n$ is an $A_n$-module: for all elements $p\in P_n$,
$$ x_i*p=x_ip, \;\;\;\;\; \der_i*p = \frac{\der p}{\der x_i},
\;\;\; i=1, \ldots , n.$$ Clearly, $P_n\simeq A_n / \sum_{i=1}^n
A_n \der_i$, $1\mapsto 1+\sum_{i=1}^n A_n \der_i$. Since
$\ggu_n\subseteq A_n$, the polynomial algebra $P_n$ is also a
$\ggu_n$-module.

Let $V$ be a vector space over $K$. A $K$-linear map $\d : V\ra V$
is called a {\em locally nilpotent map} if $V=\cup_{i\geq 1} \ker
(\d^i)$ or, equivalently, for every $v\in V$, $\d^i (v) =0$ for
all $i\gg 1$. When  $\d$ is a locally nilpotent map in $V$ we
also say that $\d$ {\em acts locally nilpotently} on $V$. Every {\em nilpotent} linear map  $\d$, that is $\d^n=0$ for some $n\geq 1$, is a locally nilpotent map but not vice versa, in general.   Let
$\CG$ be a Lie algebra. Each element $a\in \CG$ determines the
derivation  of the Lie algebra $\CG$ by the rule $\ad (a) : \CG
\ra \CG$, $b\mapsto [a,b]$, which is called the {\em inner
derivation} associated with $a$. The set $\Inn (\CG )$ of all the
inner derivations of the Lie algebra $\CG$ is a Lie subalgebra of
the Lie algebra $(\End_K(\CG ), [\cdot , \cdot ])$ where $[f,g]:=
fg-gf$. There is the short exact sequence of Lie algebras
$$ 0\ra Z(\CG ) \ra \CG\stackrel{\ad}{\ra} \Inn (\CG )\ra 0,$$
that is $\Inn (\CG ) \simeq \CG / Z(\CG )$ where $Z(\CG )$ is the {\em centre} of the Lie algebra $\CG$ and $\ad ([a,b]) = [
\ad (a) , \ad (b)]$ for all elements $a, b \in \CG$. An element $a\in \CG$ is called a {\em locally nilpotent element} (respectively, a {\em nilpotent element}) if so is the inner derivation $\ad (a)$ of the Lie algebra $\CG$. Let $J$ be a non-empty subset of $\CG$ then
$\Cen_{\CG}(J) :=\{ a\in \CG \, | \, [a,b]=0$ for all $b\in J\}$
is called the {\em centralizer} of $J$ in $\CG$. It is a Lie subalgebra of the Lie algebra $\CG$.

\begin{proposition}\label{a8Dec11}
{\rm (Proposition 2.1, \cite{Lie-Un-GEN})}
\begin{enumerate}
\item The Lie algebra $\ggu_n$ is a solvable but not nilpotent Lie
algebra. \item The finite chain of ideals (\ref{seruni}) is the
 derived  series for the Lie algebra $\ggu_n$, that is
$(\ggu_n)_{(i)}= \ggu_{n,i+1}$ for all $i\geq 0$. \item The upper
central series for the Lie algebra $\ggu_n$ stabilizers at the
first step, that is  $(\ggu_n)^{(0)}=\ggu_n$ and
$(\ggu_n)^{(i)}=\ggu_{n, 2}$ for all $i\geq 1$. \item Each element
$u\in \ggu_n$ acts locally nilpotently on the $\ggu_n$-module
$P_n$. \item All the inner derivations of the Lie algebra $\ggu_n$
are locally nilpotent derivations.  \item The centre $Z(\ggu_n)$
of the Lie algebra $\ggu_n$ is $K\der_n$. \item The Lie algebras
$\ggu_n$ where  $ n\geq 2$ are pairwise non-isomorphic.
\end{enumerate}
\end{proposition}

Proposition \ref{a8Dec11}.(5) allows us to produces many automorphisms of the Lie algebra $\ggu_n$. For every element $a\in \ggu_n$, the inner derivation $\ad (a)$ is a locally nilpotent derivation, hence $e^{\ad (a)}:=\sum_{i\geq 0} \frac{\ad (a)^i}{i!}\in G_n$. The adjoint group $\CA (\ggu_n):=\langle e^{\ad (a)}\, | \, a\in \ggu_n\rangle$ coincides with the group $\UAut_K(P_n)_n$ (Theorem \ref{8Apr12}) which is a tiny part of the group $G_n$ (Theorem \ref{28Feb12}).

{\bf The uniserial dimension}.  Let $(S, \leq )$ be a {partially
ordered set} (a {\em poset}, for short),  i.e.,   a set $S$ admits a
relation $\leq$ that satisfies three conditions: for all $a,b,c\in
S$,

(i) $a\leq a$;

(ii) $a\leq b$ and $b\leq a$ imply $a=b$;

(iii) $a\leq b$ and  $b\leq c$ imply $a\leq c$.

A poset $(S, \leq )$ is called an {\em Artinian} poset is every
non-empty subset $T$ of $S$ has a {\em minimal element}, say $t\in
T$, that is $t\leq t'$ for all $t'\in T$.  A poset $(S,\leq )$ is
a {\em well-ordered} if for all elements $a,b\in S$ either $a\leq
b$ or $b\leq a$. A bijection $f: S\ra S'$ between two posets $(S,
\leq )$ and $(S', \leq )$ is an {\em isomorphism} if $a\leq b$ in
$S$ implies $f(a)\leq f(b)$ in $S'$. Recall that the {\em ordinal
numbers} are the isomorphism classes of well-ordered Artinian
sets. The ordinal number (the isomorphism class) of a well-ordered
Artinian set $(S, \leq )$ is denoted by $\ord (S)$. The class of
all ordinal numbers is denoted by $ \mW $. The class $\mW$ is
well-ordered  by `inclusion' $\leq$ and Artinian. An associative
addition `$+$' and an associative multiplication `$\cdot$' are
defined in $\mW$ that extend the addition and multiplication of
the natural numbers.  Every non-zero natural number $n$ is
identified with $\ord (1<2<\cdots <n)$. Let $\o := \ord (\N , \leq
)$. More details on the ordinal numbers the reader can find in the book \cite{Rotman-Kneebone-Book}.

$\noindent $

{\it Definition},  \cite{Lie-Un-GEN}. Let $(S, \leq )$ be a partially ordered set. The {\em uniserial dimension} $\udim (S)$ of $S$ is the supremum of  $\ord (\CI )$ where $\CI$ runs through
all the Artinian well-ordered subsets of $S$.

$\noindent $

For a Lie algebra $\CG$, let $\CJ_0 (\CG)$ and $\CJ (\CG )$ be the
sets of all and all non-zero ideals  of the Lie algebra $\CG$,
respectively. So, $\CJ_0(\CG ) = \CJ (\CG )\cup \{ 0\}$. The sets
$\CJ_0(\CG)$ and $\CJ (\CG )$ are posets with respect to
inclusion. A Lie algebra $\CG$ is called {\em Artinian}
(respectively, {\em Noetherian}) if the poset $\CJ (\CG )$ is
Artinian (respectively, Noetherian). This means that every
descending (respectively, ascending) chains of ideals stabilizers.
A Lie algebra $\CG$ is called a {\em uniserial} Lie algebra if the
poset $\CJ (\CG )$ is a well-ordered set. This means that for any
two ideals  $\ga$ and $\gb$ of the Lie algebra $\CG$ either $\ga
\subseteq \gb$ or $\gb \subseteq \ga$.

$\noindent $

{\it Definition},  \cite{Lie-Un-GEN}. Let $\CG$ be an Artinian uniserial Lie algebra.
The ordinal number $\udim (\CG) := \ord (\CJ (\CG ))$ of the
Artinian  well-ordered set $\CJ (\CG )$ of nonzero ideals of $\CG$ is
called the {\em uniserial dimension} of the Lie algebra $\CG$. For
an arbitrary Lie algebra $\CG$, the uniserial dimension $\udim
(\CG )$ is the supremum of $\ord (\CI )$ where $\CI$ runs through
all the Artinian well-ordered sets of ideals.

$\noindent $

If $\CG$ is a Noetherian Lie algebra then $\udim (\CG ) \leq \o$.
So, the uniserial dimension is a measure of deviation from the
Noetherian condition. The concept of the uniserial dimension makes
sense for any algebras (associative, Jordan, etc.).

Let $A$ be an algebra and $M$ be its module, and let $\CJ_l(A)$ and $\CM (M)$ be the sets of all the nonzero left ideals of $A$ and of all the nonzero submodules of $M$, respectively. They are posets with respect to $\subseteq $. The {\em left uniserial dimension} of the algebra $A$ is defined as $\udim (A):= \udim (\CJ_l(A))$ and the {\em uniserial dimension} of the $A$-module $M$ is defined as $\udim (M):=\udim (\CM (M))$,  \cite{Lie-Un-GEN}.

{\bf An Artinian well-ordering on the canonical basis $\CB_n$ of
$\ggu_n$}.  Let us define an Artinian well-ordering $\leq $ on the
canonical basis $\CB_n$  for the Lie algebra $\ggu_n$ by the rule:
$\Xai
>\Xbj$ iff $i<j$ or $i=j$ and $\alpha_{n-1}=\beta_{n-1}, \ldots,
\alpha_{m+1}=\beta_{m+1}, \alpha_m>\beta_m$ for some $m$.


The next lemma is a straightforward consequence of the definition of the ordering $<$, we write $0<X_{\alpha , i}$ for all $X_{\alpha , i}\in \CB_n$.
\begin{lemma}\label{a18Dec11}
If $X_{\alpha , i}>X_{\beta , j}$ then
\begin{enumerate}
\item $X_{\alpha +\g, i}>X_{\beta +\g, j}$ for all $\g \in \N^{i-1}$,
\item $X_{\alpha -\g, i}>X_{\beta -\g, j}$ for all $\g \in \N^{i-1}$ such that $\alpha -\g \in \N^{i-1}$ and $\beta -\g \in \N^{j-1}$,
\item $[\der_k, X_{\alpha , i}]>[\der_k, X_{\beta , j}]$ for all $k=1, \ldots , i-1$ such that $\alpha_k\neq 0$, and
    \item $[X_{\g ,k}, X_{\alpha , i}]>[X_{\g ,k}, X_{\beta , j}]$ for all $X_{\g , k}>X_{\alpha , i}$ such that
           $[X_{\g ,k}, X_{\alpha , i}]\neq 0$, i.e.,  $\alpha_k\neq 0$.
\end{enumerate}
\end{lemma}

Let $\O_n$ be the set of indices $\{ (\alpha , i)\}$ that
parameterizes the canonical basis $\{ \Xai \}$ of the Lie algebra
$\ggu_n$. The set $(\O_n, \leq )$ is an Artinian well-ordered set,
where $(\alpha , i)\geq (\beta , j)$ iff $\Xai \geq \Xbj$, which
is isomorphic to the Artinian well-ordered set $(\CB_n, \leq )$
via $(\alpha , i)\mapsto \Xai$. We identify the posets $(\O_n,
\leq )$ and $(\CB_n , \leq )$ via this isomorphism.  It is obvious
that 
\begin{equation}\label{ordn}
\ord (\CB_n)=\ord (\O_n ) = \o^{n-1}+\o^{n-2}+\cdots +\o +1,
\end{equation}
$\O_2\subset \O_3\subset \cdots $ and $\CB_2\subset \CB_3\subset
\cdots $.   Let $[1, \ord (\O_n)]:=\{ \l \in \mW \, | \, 1\leq \l
\leq \ord (\O_n)\}$. By (\ref{Xai1}), if $[\Xai, \Xbj ]\neq 0$
then 
\begin{equation}\label{XaXj}
[\Xai , \Xbj ] <\min \{ \Xai , \Xbj \}.
\end{equation}

{\bf A classification of ideals of the Lie algebra $\ggu_n$}.
 By (\ref{XaXj}), the  map 
\begin{equation}\label{Jlid}
\rho_n: [1, \ord (\O_n)]\ra \CJ (\ggu_n), \;\; \l  \mapsto I_\l :=I_{\l, n} :=
\bigoplus_{ (\alpha , i)\leq \l } K\Xai ,
\end{equation}
is a monomorphism  of posets ($\rho_n$  is an order-preserving
injection).

\begin{theorem}\label{10Dec11}
{\rm (Theorem 3.3, \cite{Lie-Un-GEN})}
\begin{enumerate}
\item The map (\ref{Jlid}) is a bijection. \item The Lie algebra
$\ggu_n$ is a uniserial, Artinian but not Noetherian Lie algebra
and its uniserial dimension is equal to $\udim (\ggu_n) = \ord
(\O_n) = \o^{n-1}+\o^{n-2}+\cdots +\o +1$.
\end{enumerate}
\end{theorem}
An ideal $\ga$ of a Lie algebra $\CG$ is called {\em proper}
(respectively, {\em co-finite}) if $\ga \neq 0, \CG$
(respectively, $\dim_K(\CG / \ga ) <\infty $).
An ideal $I$  of a Lie algebra $\CG$ is called a {\em
characteristic ideal} if it is invariant under all the
automorphisms of the Lie algebra $\CG$, that is $\s (I)=I$ for all
$\s \in \Aut_K(\CG )$. It is obvious that an ideal $I$ is a
characteristic ideal iff $\s (I) \subseteq I$ for all $\s \in
\Aut_K(\CG )$.

\begin{corollary}\label{c10Dec11}
{\rm (Corollary 3.7, \cite{Lie-Un-GEN})}
All the ideals of the Lie algebra $\ggu_n$ are characteristic
ideals.
\end{corollary}

 Each
non-zero element $u$ of $\ggu_n$ is a finite linear combination
$$ u = \l\Xai +\mu \Xbj +\cdots +\nu X_{\s , k}= \l \Xai +\cdots $$
where $\l , \mu , \ldots , \nu \in K^*$ and $\Xai > \Xbj >\cdots
>X_{\s , k}$. The elements $\l \Xai$ and $\l \in K^*$ are called
the {\em leading term} and the {\em leading coefficient} of $u$
respectively, and the ordinal number  $\ord (\Xai ) = \ord (\alpha
, i) \in [ 1, \ord (\O_n)]$, which is, by definition, the ordinal number that represents the Artinian well ordered set $\{ (\beta, j) \in \O_n \, | \, (\beta , j)\leq (\alpha , i)\}$,  is called the {\em ordinal degree} of
$u$ denoted by $\ord (u)$ (we hope that this notation will not
lead to confusion). For all non-zero elements $u,v\in \ggu_n$ and
$\l \in K^*$,

(i) $\ord (u+v) \leq \max \{ \ord (u) , \ord (v)\}$ provided
$u+v\neq 0$;

(ii) $\ord (\l u) = \ord (u)$;

(iii) $\ord ([u,v])<  \min \{ \ord (u), \ord (v)\} $ provided
$[u,v]\neq 0$;


\begin{corollary}\label{d10Dec11}
{\rm (Corollary 3.8, \cite{Lie-Un-GEN})}
For all nonzero elements $u\in \ggu_n$ and all automorphisms $\s $
of the Lie algebra $\ggu_n$, $\ord (\s (u)) = \ord (u)$.
\end{corollary}




\section{The structure of the group of automorphisms of the  Lie algebra $\ggu_n$}\label{GALAUN}

In this section, several important subgroups of the group $G_n$ are introduced and studied. It is proved that the group $G_n$ is an iterated semi-direct product and an exact product of some of them (Proposition \ref{17Dec11}, Theorem \ref{26Jan12}).

{\bf The groups $\mT^n$ and $\CU_n$}.
Let $G_n:=\Aut_K(\ggu_n)$ be the group of automorphisms of the Lie algebra $\ggu_n$. The Lie algebra $\ggu_n$ is uniserial, so 
\begin{equation}\label{sIl}
\s (I_\l ) = I_\l \;\; {\rm for \; all}\;\; \l \in [1, \ord (\O_n)],
\end{equation}
by Theorem \ref{10Dec11}. Moreover, by (\ref{Jlid}),
\begin{equation}\label{sXl}
\s (X_{\alpha , i} ) =\l_{\alpha , i} X_{\alpha , i}+\cdots  \;\; {\rm for \; all}\;\; X_{\alpha , i}\in \CB_n
\end{equation}
for some scalar $\l_{\alpha , i}=\l_{\alpha , i}(\s )\in K^*$ where the three dots mean smaller terms, i.e.,  an element of $\sum_{(\beta , j) <(\alpha , i)}K X_{\beta , j}$. It follows that
$$\CU_n:= \{ \s \in G_n \, | \, \s (X_{\alpha , i})  = X_{\alpha , i}+\cdots\;\; {\rm for \; all}\;\; X_{\alpha, i}\in \CB_n\}$$ is a normal subgroup of the group $G_n$. The {\em algebraic $n$-dimensional torus} $\mT^n$ is a subgroup of the group $\Aut_K(A_n)$ of automorphisms of the Weyl algebra $A_n$,
$$ \mT^n := \{ t_\l : x_i\mapsto \l_ix_i, \; \der_i\mapsto \l_i^{-1}\der_i, 1\leq i\leq n\, | \, \l = (\l_i) \in K^{*n}\} \simeq K^{*n}, $$ that preserves the Lie algebra $\ggu_n$. The group $\mT^n$ can be seen as a subgroup of $G_n$,
$$ \mT^n := \{ t_\l : X_{\alpha, i}\mapsto \l^\alpha \l_i^{-1}X_{\alpha, i}\; {\rm for \; all } \; X_{\alpha , i}\in \CB_n \, | \, \l = (\l_i) \in K^{*n}\} \simeq K^{*n}$$
where $\l^\alpha := \prod \l_i^{\alpha_i}$.
\begin{proposition}\label{17Dec11}
\begin{enumerate}
\item The group $\CU_n$ is a normal subgroup of  $G_n$ and $\CU_n=\{ \s \in G_n\, | \, \s (\der_i) = \der_i+\cdots $ for $i=1,\ldots , n\}$.
\item $G_n=\mT^n\ltimes \CU_n$ (the group $G_n$ is the semidirect product of $\mT_n$ and $\CU_n$).
\end{enumerate}
\end{proposition}

{\it Proof}. 1. We have already seen that $\CU_n$ is a normal subgroup of $G_n$. It remains to show that the equality holds. Let $R$ be the RHS of the equality. Then $\CU_n\subseteq R$. It remains to show that $\CU_n\supseteq R$, i.e.,  $\s \in R$ implies $\s \in \CU_n$. We have to show that $\s (X_{\alpha , i})=X_{\alpha , i}+\cdots $ for all $X_{\alpha , i}\in \CB_n$.
 We use induction on $\l := \ord (X_{\alpha , i}) = \ord ((\alpha , i))\in [ 1, \ord (\O_n)]$. The initial case $\l =1$ is obvious as $X_{\alpha ,i}=\der_n$ and  $\s (\der_n ) =\der_n$ since $\s \in R$.

Let $\l >1$, and we assume that the result is true for all $\l'<\l $. If $X_{\alpha , i}=\der_j$ for some $j$ then it is nothing to prove. So, let $X_{\alpha ,  i}\in \CB_n\backslash \{ \der_1, \ldots , \der_n\}$, i.e.,  $\alpha \in \N^{i-1}\backslash \{ 0\}$. Let $j=\max \{ k\, | \, \alpha_k\neq 0\}$. Then, by the very definition of the  ordering $<$ on $\CB_n$ (or use Lemma \ref{a18Dec11}.(3)),
$$ [\der_j, X_{\alpha , i}+\cdots ] = \alpha_jX_{\alpha -e_j, i}+\cdots \;\; {\rm and}\;\; [\oplus_{k=j+1}^iP_{k-1}\der_k , X_{\alpha , i}]=0.$$
Then applying the automorphism $\s$ to the identity $[\der_j, X_{\alpha , i}]= \alpha_j X_{\alpha -e_j, i}$ and using the fact that $\s (\der_j) = \der_j+u+v$ for some elements $u\in \oplus_{k=j+1}^iP_{k-1}\der_k$ and $v\in \oplus_{k=i+1}^nP_{k-1}\der_k$, we have
\begin{eqnarray*}
 \alpha_jX_{\alpha -e_j, i}+\cdots & =& \s (\alpha_j X_{\alpha -e_j, i})=[\s (\der_j), \s (X_{\alpha , i})]= [ \der_j+u+v, \l_{\alpha ,i}X_{\alpha , i}+\cdots ]     \\
 & =&\l_{\alpha , i}\alpha_jX_{\alpha -e_j, i}+[v,\l_{\alpha , i}X_{\alpha , i}]+[\der_j+u+v, \cdots ]= \l_{\alpha , i}\alpha_j X_{\alpha -e_j, i}+\cdots .
\end{eqnarray*}
Hence, $\l_{\alpha , i}=1$ since $\alpha_j\neq 0$.

2. By the very definition of the groups $\mT^n$ and $\CU_n$, $\mT^n \cap \CU_n = \{ e\}$. As $\CU_n$ is a normal subgroup of $G_n$ it  suffices to show that $G_n=\mT^n \CU_n$, i.e.,  every automorphism $\s \in G_n$ is a product $t\tau $ for some elements $t\in \mT^n$ and $\tau \in \CU_n$. By (\ref{sXl}), $\s (\der_i) = \l_i\der_i+\cdots $ for all $i=1, \ldots , n$ and for some $\l = (\l_1, \ldots , \l_n) \in K^{*n}$. Then $t_\l \s (\der_i) = \der_i+\cdots$ for all $i=1, \ldots , n$ where $t_\l \in \mT^n$. By statement 1, $\tau := t_\l \s \in \CU_n$, hence $\s = t_\l^{-1}\tau$, as required.  $\Box $


{\bf The group $\TAut_K(P_n)$ of triangular automorphisms of the polynomial algebra  $P_n$}. Let $\Aut_K(P_n)$ be the group of $K$-algebra automorphisms of the polynomial algebra $P_n$. Every automorphism $\s \in \Aut_K(P_n)$ is uniquely determined by the polynomials
$$x_1':=\s (x_1), \ldots , x_n':= \s (x_n).$$ The inclusions of the polynomial algebras $P_1\subset P_2\subset \cdots$ yield the natural inclusions of their  groups of automorphisms $\Aut_K(P_1)\subset \Aut_K(P_2)\subset \cdots$ where an automorphism $\s \in \Aut_K(P_n)$ is extended to the automorphism of the polynomial algebra $P_{n+1}$ by the rule $\s (x_{n+1}) = x_{n+1}$.

 The {\em group of triangular automorphisms} $\TAut_K(P_n)$ of the polynomial algebra $P_n$ consists of all the  automorphisms of $P_n$   of the following type:
\begin{equation}\label{sa1n}
\s (x_i)=\l_ix_i+a_i,\;\; i=1, \ldots , n,
\end{equation}
where $a_i\in P_{i-1}$ and $\l_i\in K^*$ for $i=1, \ldots , n$.
The automorphism $\s$ is uniquely determined by the elements $a_i$
and $\l_i$, and we write $\s = [a_1, \ldots , a_n ; \l_1, \ldots ,
\l_n]$. There are two distinct subgroups in $\TAut_K(P_n)$: the
{\em algebraic $n$-dimensional torus} $\mT^n$ (where $a_1=\cdots =
a_n=0$)  and the {\em group $\UAut_K(P_n)$ of triangular
polynomial automorphisms} (where $\l_1=\cdots = \l_n=1$).
Moreover, $\UAut_K(P_n)$ is a normal subgroup of $\TAut_K(P_n)$
and
\begin{equation}\label{TAut=TU}
\TAut_K(P_n) = \mT^n\ltimes \UAut_K(P_n).
\end{equation}

{\bf The group $\UAut_K(P_n)$ of triangular automorphisms of  the
polynomial algebra $P_n$}. Every element $[a_1, \ldots , a_n]:=
[a_1, \ldots , a_n; 1, \ldots , 1]\in \UAut_K(P_n)$ is uniquely
determined by the polynomials $a_i\in P_{i-1}$, $i=1, \ldots , n$.
\begin{proposition}\label{a28Jan12}
The exponential map $ \ggu_n\ra \UAut_K(P_n)$, $ \d \mapsto e^\d:=\sum_{i\geq 0} \frac{\d^i}{i!}$, is a bijection with the inverse map $ \s\mapsto \ln (\s) := \ln (1-(1-\s)):=-\sum_{i\geq 1} \frac{(1-\s )^i}{i}$.
\end{proposition}

{\it Proof}. By Proposition \ref{a8Dec11}.(4), the exponential map is well-defined. Let us show that for every automorphism $\s = [a_1, \ldots , a_n]\in \UAut_K(P_n)$ there is the unique derivation $\d = \sum_{i=1}^n b_i\der_i\in \ggu_n$ (where $b_i\in P_{i-1}$ for $i=1, \ldots , n$) such that $ \s = e^\d$. Consider the system of equations where $\{ b_i \}$ are unknown polynomials such that $\s (x_i) = e^\d (x_i)$, $i=1, \ldots , n$, that is
$x_i +a_i = x_i +(1-\der ) (b_i) $ where $\der := -\sum_{i\geq 1} \frac{\d^i}{(i+1)!}$ is a locally nilpotent map on $P_n$. Then
$$ b_i = (1-\der )^{-1} (a_i) = (\sum_{j\geq 0} \der^j ) (a_i)\in P_{i-1}.$$
For each $\s \in \UAut_K(P_n)$, the map $ 1-\s : P_n\ra P_n$ is a locally nilpotent map. So, the map  $\ln (\s )= -\sum_{i\geq 1} \frac{(1-\s )^i}{i}$ makes sense. The rest is obvious. $\Box $


For every element $[a_1, \ldots , a_n] \in \UAut_K(P_n)$,
\begin{equation}\label{a1e}
[a_1, \ldots , a_n]=e^{a_n\der_n} e^{a_{n-1}\der_{n-1}}\cdots e^{a_1\der_1}.
\end{equation}
For each natural number $i=1, \ldots , n$, the map $P_{i-1}\ra e^{P_{i-1}\der_i}$, $p_i\mapsto e^{p_i\der_i}$, is an isomorphism of abelian groups. The group $\UAut_K(P_n)$ is an iterated semi-direct product of its subgroups $e^{P_{i-1}\der_i}$, $i=1, \ldots , n$,
\begin{equation}\label{a1ex1}
\UAut_K(P_n) = e^{P_{n-1}\der_n}\rtimes e^{P_{n-2}\der_{n-1}}\rtimes \cdots \rtimes e^{P_0\der_1}.
\end{equation}
The set of all $K$-derivations $\Der_K(P_n)=\oplus_{i=1}^n P_n\der_n$ of the polynomial algebra $P_n$ is a Lie subalgebra of the Lie Weyl algebra $(A_n, [\cdot , \cdot ])$, and $\ggu_n$ is a Lie subalgebra of $\Der_K(P_n)$. Each automorphism $\s \in \Aut_K(P_n)$ induces an  automorphism of the Lie algebra $\Der_K(P_n)$ (the change of variable) by the rule $\d \mapsto \s \d \s^{-1}$ where $\d \in \Der_K(P_n)$. In particular, $\s (\frac{\der}{\der x_i}) \s^{-1} = \frac{\der}{\der x_i'}$ where $x_i':= \s (x_i)$.  Moreover,
\begin{equation}\label{derjp}
\der_i':=\frac{\der}{\der x_i'}=\sum_{j=1}^n\s(\frac{\der\s^{-1}(x_j)}{\der x_i})\der_j.
\end{equation}
Let $\s \in \TAut_K(P_n)$ be as in (\ref{sa1n}). Then the automorphism $\s^{-1}\in \TAut_K(P_n)$ and so $\s^{-1} = (b_1, \ldots , b_n; \l_1^{-1}, \ldots , \l_n^{-1})$ for some $b_i\in P_{i-1}$, $i=1, \ldots , n$. For the automorphism $\s$ the equality (\ref{derjp}) takes form
\begin{equation}\label{derjp1}
\der_i'=\frac{\der}{\der x_i'}=\l_i^{-1}\der_i+\sum_{j=i+1}^n\s(\frac{\der\s^{-1}(b_j)}{\der x_i})\der_j\in K^*\der_i+\sum_{j=i+1}^n P_{j-1}\der_j.
\end{equation}
{\bf The groups $\Sh_n$, $\CT_n$ and $\UAut_K(P_n)_n$}. Important subgroups of $\TAut_K(P_n)$ are the {\em shift group}
\begin{equation}\label{Shn}
\Sh_n:= \{ \s : x_1\mapsto x_1+\mu_1, \ldots , x_n\mapsto x_n+\mu_n \, | \; \mu_i\in K, i=1, \ldots , n\}\simeq K^n
\end{equation}
and the  group
\begin{equation}\label{CTndef}
\CT_n:=\{ \s = [0, a_2, \ldots , a_n; 1, \ldots , 1] \, | \, a_2(0)=a_3(0,0)=\cdots = a_n(0, \ldots , 0)=0\}.
\end{equation}
So, an automorphism $\s = [a_1, a_2, \ldots , a_n; \l_1, \ldots , \l_n ]\in \TAut_K(P_n)$ belongs to the group $\CT_n$ iff $\l_1=\cdots = \l_n=1$ and all the constant terms of the polynomials $a_i$ are equal to zero. Notice that $\CT_n\subseteq \UAut_K(P_n)\subseteq  \TAut_K(P_n)$. For each $i\geq 1$, let $\gm_i:= \sum_{j=1}^ix_jP_i$, the maximal ideal of the polynomial algebra $P_i$ generated by the elements $x_1, \ldots , x_i$. Set $\gm_0:=0$. So, an automorphism $\s \in \UAut_K(P_n)$ belongs to the group  $\CT_n$ iff $\s (\gm_i) = \gm_i$ for all $i=1, \ldots , n$.  Let
\begin{eqnarray*}
 \TAut_K(P_n)_n &:=&\{ \s \in \TAut_K(P_n)\, | \, \s (x_n) = \l_n x_n+a_n\;\; {\rm where}\;\;\l_n\in K^*, \;  a_n\in \gm_{n-1}\},  \\
 \UAut_K(P_n)_n &:=&\{ \s \in \UAut_K(P_n)\, | \, \s (x_n) = x_n+a_n\;\; {\rm where}\;\; a_n\in \gm_{n-1}\},  \\
 \sh_n&:=&\Sh_n\cap \Fix_{\TAut_K(P_n)}(x_1, \ldots , x_{n-1})=\{ \s \in \Sh_n \, | \, \s(x_i) = x_i, \; i=1, \ldots , n-1; \\
  & &  \s (x_n) = x_n+\l , \l \in K\}
 = e^{K\der_n}:= \{ e^{\l \der_n}\, | \, \l \in K\}\simeq (K, +).
\end{eqnarray*}
The subsets $\TAut_K(P_n)_n$ and $\UAut_K(P_n)_n$ of the group $\TAut_K(P_n)$ are not subgroups, but $\sh_n$ is a subgroup.
Clearly, $\TAut_K(P_n) = \TAut_K(P_n)_n\sh_n$, $\TAut_K(P_n)_n\cap \sh_n=\{ e\}$ and the group $\sh_n$ is a normal subgroup  of the group  $\TAut_K(P_n)$.
 Therefore, the sets $\TAut_K(P_n)_n$ and $\UAut_K(P_n)_n$ can be identified with the factor groups $\TAut_K(P_n)/\sh_n$  and $\UAut_K(P_n)/\sh_n$, respectively, and as a result they have the group structure. Under these identifications we can write
\begin{equation}\label{derjp2}
\TAut_K(P_n)_n = \TAut_K(P_n)/\sh_n,
\end{equation}
\begin{equation}\label{derjp3}
\UAut_K(P_n)_n = \UAut_K(P_n)/\sh_n.
\end{equation}
\begin{proposition}\label{a27Jan12}

\begin{enumerate}
\item $\s\ggu_n\s^{-1} = \ggu_n$ for all $\s \in \TAut_K(P_n)$.
\item The map $\o : \TAut_K(P_n) \ra G_n$, $\s \mapsto (\o_\s : u\mapsto \s u\s^{-1})$, (where $u\in \ggu_n$)  is a group homomorphism with $\ker (\o ) = \sh_n$.
\item The map $\o : \TAut_K(P_n)_n \ra G_n$, $\s \mapsto \o_\s$, is a group monomorphism.

\end{enumerate}
\end{proposition}

{\it Proof}. 1. Since $\TAut_K(P_n)$ is a group, to prove statement 1 it suffices to show that $\s \ggu_n\s^{-1} \subseteq \ggu_n$ for all elements $\s \in \TAut_K(P_n)$. Since $\ggu_n = \sum_{i=1}^nP_{i-1}\der_i$, it suffices to show that $\s P_{i-1} \der_i\s^{-1} \subseteq \ggu_n$ for all elements $\s \in \TAut_k(P_n)$ and $i=1, \ldots , n$. This follows from (\ref{derjp1}) and the fact that $ \s (P_{i-1}) = P_{i-1}$:
$$ \s P_{i-1}\der_i\s^{-1} = \s (P_{i-1}) \s \der_i \s^{-1} \subseteq P_{i-1} (K^*\der_i+\sum_{j=i+1}^nP_{j-1}\der_j) \subseteq \sum_{j=i}^nP_{j-1}\der_j\subseteq \ggu_n.$$

2. The map $\o$ is  a group homomorphism. By (\ref{derjp}), $\sh_n\subseteq \ker (\o )$. Let $\s\in \ker (\o )$. It remains to show that $\s \in \sh_n$. By (\ref{derjp1}), $\der_i'=\der_i$ for all $i=1, \ldots , n$. Hence $\l_1=\cdots = \l_n=1$ and $\frac{\der \s^{-1}(b_j)}{\der x_i}=0$ for all $i\leq j$. The second set of the  conditions means that the elements $\s^{-1} (b_j) \in P_{j-1}$ are scalars. Summarizing, $\s^{-1} (x_i) = x_i+b_i$ where all $b_i\in K$. For all $i=2, \ldots , n$,
$$ x_{i-1}\der_i = \s^{-1} ( x_{i-1}\der_i)\s  = \s^{-1} (x_i) \der_i = x_{i-1}\der_i +b_{i-1}\der_i.$$ Hence, $ b_1=\cdots  = b_{n-1}=0$. Therefore, $\s \in \sh_n$, as required.

3. Statement 3 follows from statement 2 and (\ref{derjp2}).  $\Box $


By Proposition \ref{a27Jan12}.(3), we identify the group $\TAut_K(P_n)_n=\TAut_K(P_n)/\sh_n$ with its image in the group $G_n$, i.e.,  $\TAut_K(P_n)_n\subseteq G_n$. We identify the subgroup $\Sh_{n-1}$ of $\Aut_K(P_{n-1})$ with the following subgroup of $\Aut_K(P_n)$, 
\begin{equation}\label{Shn1}
\Sh_{n-1} = \{ \s \in \Sh_n \, | \, \s (x_n) = x_n \}.
\end{equation}
$$ \Sh_n=\sh_1\times \cdots \times \sh_n = e^{K\der_1}\times\cdots \times e^{K\der_n}=e^{\sum_{i=1}^nK\der_i}.$$

We say that a group $G$ is the {\em exact product} of its two (or several) subgroups $G_1$ and $G_2$ and write $G= G_1\times_{ex}G_2$ if every element $g\in G$ is a unique product $g=g_1g_2$ where $g_1\in G_1$ and $g_2\in G_2$. Using the bijection $G\ra G$, $g\mapsto g^{-1}$ and the fact that $(gh)^{-1} = h^{-1}g^{-1}$, we see that $G= G_1\times_{ex}G_2$ iff $G= G_2\times_{ex}G_1$. The next theorem describes the group $G_n$.

\begin{proposition}\label{c27Jan12}

\begin{enumerate}
\item $\UAut_K(P_n)=\CT_n\times_{ex}\Sh_n$.
\item $\TAut_K(P_n)=\mT^n\ltimes ( \CT_n\times_{ex}\Sh_n)$.
\item $\TAut_K(P_n)_n= \mT^n\ltimes ( \CT_n\times_{ex}\Sh_{n-1})=\mT^n\ltimes \UAut_K(P_n)_n$.
\item $\UAut_K(P_n)_n= \CT_n\times_{ex}\Sh_{n-1}$.
\end{enumerate}
\end{proposition}

{\it Proof}. 1.  It is obvious that $ \Sh_n\cap \CT_n = \{ e\}$. To finish the proof of statement 1 it suffices to show that any automorphism $ \s = [a_1, a_2, \ldots , a_n; 1, \ldots , 1 ]\in \UAut_K(P_n)$ is a product $\tau s$ for some elements $ \tau \in \CT_n$ and $s\in \Sh_n$. Let $a_i = b_i+\mu_i$ where $ \mu_i$ is the constant term of the polynomial $a_i$. Then  $\tau =  [0, b_2, \ldots , b_n; 1, \ldots , 1 ]\in \CT_n$,  $s= [\mu_1, \ldots , \mu_n; 1, \ldots , 1 ]\in \Sh_n$ and $ \s = \tau s$.

2. Statement 2 follows from statement 1 and (\ref{TAut=TU}).

3. Notice that $\Sh_n= \Sh_{n-1}\times \sh_n$. Statement 3 follows from statement 1 and (\ref{derjp2}).

4. Statement 4 follows from statement 1 and (\ref{derjp3}). $\Box$


The following lemma gives a characterization of the shift group $\Sh_n$ via its action on the partial derivatives.

\begin{lemma}\label{b27Jan12}
$\Fix_{\Aut_K(P_n)} (\der_1, \ldots , \der_n) = \Sh_n$ and $\Fix_{\UAut_K(P_n)_n} (\der_1, \ldots , \der_n) = \Sh_{n-1}$.
\end{lemma}

{\it Proof}. For an automorphism $\s \in \Aut_K(P_n)$, let $\der':= ( \der_1', \ldots , \der_n')^T$ where $\der_i':= \s \der_i\s^{-1}$ for $i=1, \ldots , n$,  $\der:= ( \der_1, \ldots , \der_n)^T$ (where `T' stands for the transposition)  and $A=(A_{ij})$ be the $n\times n$ matrix where $A_{ij}:= \s (\frac{\der \s^{-1}(x_j)}{\der x_i})$. The equalities (\ref{derjp}) can be written in the matrix form as $ \der' = A\der$. Then $\s \in \Fix_{\Aut_K(P_n)} (\der_1, \ldots , \der_n)$ iff $A$ is the identity matrix iff $\s \in  \Sh_n$. The second equality follows from the first. $\Box$


The next theorem is a key result which is used in several proofs later.
\begin{theorem}\label{B27Jan12}
 Let  $\der_1', \ldots , \der_n'\in \Der_K(P_n)$ be commuting derivations such that $\der_i':= \mu_i \der_i +\sum_{j=i+1}^n a_{ij} \der_j$ for $i=1, \ldots , n$ where $\mu_i\in K^*$ and $a_{ij} \in P_{j-1}$. Then
\begin{enumerate}
\item there exists an automorphism $\s \in \TAut_K(P_n)$ such that $\der_i'= \s \der _i \s^{-1}$ for $i=1, \ldots , n$. If $\s'$ is another such an automorphism then $\s' = \s s$ for some $\s\in \Sh_n$, and vice versa.
\item There is the unique automorphism $\s \in \mT^n\ltimes \CT_n$ such that $\der_i'= \s \der_i \s^{-1} $ for $i=1, \ldots , n$. The automorphism $\s$ is defined as follows $\s (x_i) = x_i'$ for $i=1, \ldots , n$ where the elements $x_i'$ are defined recursively as follows:
\begin{equation}\label{xip}
x_1':=\mu_1^{-1}x_1, \;\; x_i':=\phi_{i-1}\phi_{i-2}\cdots \phi_1 (\mu_i^{-1}x_i) , \;\; i=2, \ldots , n,
\end{equation}
\begin{equation}\label{xip1}
\phi_i:= \sum_{k\geq 0}(-1)^k \frac{x_i'^k}{k!}\der_i'^k, \;\; i=1, \ldots , n-1.
\end{equation}
\end{enumerate}
\end{theorem}

{\it Proof}. 1. The first part of statement 1 (concerning the existence of $\s$)  follows from statement 2. If $\der_i'= \s \der_i\s^{-1} =  \s' \der_i\s'^{-1}$ for $i=1, \ldots , n$ then $s:= \s^{-1} \s'\in \Sh_n$, by Lemma \ref{b27Jan12}. Hence, $\s'= \s s$, and vice versa (trivial).

2.  We deduce  statement 2 from two claims below. The uniqueness of $\s$ follows from the second part of statement 1 (which has already been proved above) and
the fact that $\mT^n\ltimes \CT_n\cap \Sh_n = \{ e\}$.

{\em Claim 1: $\der_1', \ldots , \der_n'$ are commuting, locally nilpotent derivations of the polynomial algebra $P_n$ such that $\der_i'(x_j') = \d_{ij}$ (the Kronecker delta) and $\cap_{i=1}^n \ker_{P_n} (\der_i') = K$.}

It follows from Claim 1 and  (Theorem 2.2, \cite{Bav-inform}) that the $K$-algebra homomorphism
$$\s : P_n\ra P_n, \;\; x_i\mapsto x_i',\;\;  i=1, \ldots , n,$$ is an automorphism. Then  (which is obvious)  $\der_i'= \s \der_i \s^{-1}$ for $i=1, \ldots , n$. We finish the proof of statement 2, by Claim 2.

{\em Claim 2:} $\s \in \mT^n\ltimes \CT_n$.

So, it remains to prove Claims 1 and 2.
 By the very definition of the derivations $\der_i'$, the following statements are obvious:

(i) $\der_1', \ldots , \der_n'\in \ggu_n$;

(ii) $\der_i'(P_j) \subseteq P_j$ for all $i,j=1, \ldots , n$. Moreover, $\der_i'(P_j) =0$ for $i> j$;

(iii) $\der_i'(\mu_i^{-1}x_i) = 1$ for $i=1, \ldots , n$;

(iv) $\der_1', \ldots , \der_n'$ are locally nilpotent derivations of the polynomial algebra $P_n$;

(v) $\Der_K(P_n) = \oplus_{i=1}^n P_n \der_i'$; and

(vi) $\cap_{i=1}^n \ker_{P_n}(\der_i') = K$, by (v).

In view of (iv) and (vi), to finish the proof of Claim 1, we have to show that $\der_i'( x_j') = \d_{ij}$.  To prove Claim 2, it suffices to show that $x_j'=\mu_j^{-1} x_j +a_j$ where $a_j\in \gm_{j-1}:= \sum_{k=1}^{j-1} x_k P_{j-1}$ for $j=1, \ldots , n$.
 To prove both  statements we use induction on $j$. The initial case $j=1$ follows from (ii) and (iii) and the fact that $x_1'= \mu_1^{-1}x_1$. Suppose that $j\geq 2$ and the result  holds for all $j'<j$. By induction, $x_s'=\mu_s^{-1} x_s+a_s\in K^*x_s+\gm_{s-1} \subseteq \gm_s$ for $s=1, \ldots , j-1$. Then, using repeatedly (ii) and (\ref{xip1}), we see that $x_j'=\mu_j^{-1} x_j +a_j$ for some $a_j\in \gm_{j-1}$.  Then $\der_{j+1}'(x_j')=\cdots = \der_n'(x_j') =0$, by (ii). By  (Theorem 2.2, \cite{Bav-inform}), $\der_1'(x_j') =\cdots = \der_{j-1}'(x_j') =0$. For each $i=1, \ldots , j-1$, $\der_j'\phi_i = \phi_i\der_j'$ since $\der_j'(x_i')=0$ for all $i=1, \ldots , j-1$ (the last set of equalities follows from the set of equalities $\der_j'(x_i)=0$ for all $i=1, \ldots, j-1$, and the fact that $x_s'=\mu_s^{-1}x_s+a_s$ for $s=1, \ldots , j-1$). Therefore (by (iii)),
 $$ \der_j'(x_j') = \der_j'\phi_{j-1}\cdots \phi_1(\mu_j^{-1}x_j) = \phi_{j-1}\cdots \phi_1\der_j'(\mu_j^{-1}x_j) =\phi_{j-1}\cdots \phi_1(1) =1.$$The proof of the inductive step is complete.  $\Box $


The next theorem states that the triangular polynomial automorphisms are {\em precisely} the polynomial automorphisms that respect the Lie algebra $\ggu_n$.

\begin{theorem}\label{A27Jan12}
Let $\s\in \Aut_K(P_n)$. The following statements are equivalent.
\begin{enumerate}
\item $\s \in \TAut_K(P_n)$.
\item $\s \ggu_n\s^{-1}= \ggu_n$
\item For all $i=1, \ldots , n$, $\s \der_i \s^{-1} = \mu_i \der_i +\sum_{j=i+1}^n a_{ij}\der_j$ where $\mu_i\in K^*$ and $a_{ij} \in P_{j-1}$.
\end{enumerate}
\end{theorem}

{\it Proof}. 1. $(1\Rightarrow 2)$ Proposition \ref{a27Jan12}.(1).

$(2\Rightarrow 3)$ Suppose that $\s \ggu_n \s^{-1} = \ggu_n$. Then $\s^{-1}\ggu_n \s = \ggu_n$. The map $\o_\s : \ggu_n \ra \ggu_n$, $\d \mapsto \s \d \s^{-1}$, is a Lie algebra automorphism with $\o_{\s^{-1}}$ as its inverse. Statement 3 follows from (\ref{sIl}) and the definition of the ordering $<$ on $\ggu_n$.

$(3\Rightarrow 1)$ By Theorem \ref{B27Jan12}, there is an automorphism $\tau \in \TAut_K(P_n)$ such that $\s \der_i \s^{-1} = \tau \der_i \tau^{-1}$ for all $i=1, \ldots , n$. Hence $\tau^{-1} \s \in \Fix_{\Aut_K(P_n)} (\der_1, \ldots , \der_n) =\Sh_n\subseteq \TAut_K(P_n)$ (Lemma \ref{b27Jan12}). Then $\s \in \TAut_K(P_n)$. $\Box $


{\bf The group $G_n$ is an exact product of its three subgroups}. The group $\CU_n$ contains the subgroup
\begin{equation}\label{Fixdd}
\CF_n := \Fix_{G_n} (\der_1, \ldots , \der_n) := \{ \s \in G_n \, | \, \s (\der_1) = \der_1, \ldots , \s (\der_n) =\der_n \}.
\end{equation}
This is the most important subgroup of $G_n$ as Theorem \ref{26Jan12}.(2) shows. The next theorem represents the group $G_n$ as an exact product of its three subgroups.

\begin{theorem}\label{26Jan12}

\begin{enumerate}
\item $G_n=\TAut_K(P_n)_n \CF_n = \CF_n \TAut_K(P_n)_n$ and $\TAut_K(P_n)_n \cap \CF_n = \Sh_{n-1}$.
\item $G_n=\mT^n\ltimes (\CT_n\times_{ex} \CF_n )=\mT^n\ltimes ( \CF_n \times_{ex} \CT_n)$.
    \item $\CU_n=\CT_n\times_{ex} \CF_n$.
\end{enumerate}
\end{theorem}

{\it Proof}. 2. Let $g\in G_n$. The elements $\der_1, \ldots , \der_n$ of the Lie algebra $\ggu_n$ commute then so do the elements $\der_1':= g(\der_1) , \ldots , \der_n':=g(\der_n)$ of $\ggu_n$. By (\ref{sXl}), the elements $\der_1', \ldots , \der_n'$ satisfy the assumptions of Theorem \ref{B27Jan12}, hence there is the {\em unique} automorphism $\s \in \mT^n\ltimes \CT_n$ such that $\der_i'= \s \der_i \s^{-1}$ for $i=1, \ldots , n$. Recall that we identified the group $\mT^n\ltimes \CT_n$ with its image in the group $G_n$ (Proposition \ref{a27Jan12}.(3)), i.e.,  the automorphism $\s$ is identified with the automorphism $\o_\s : \ggu_n\ra \ggu_n$, $ u\mapsto \s u \s^{-1}$. Then $\o_{\s^{-1}}g(\der_i) = \der_i $ for all $i=1, \ldots , n$. Hence, $\o_{\s^{-1}}g\in \CF_n$, and statement 2 follows.

1. Proposition \ref{a27Jan12}.(3), the inclusion $\mT^n \ltimes \CT_n \subseteq \TAut_K(P_n)_n$ (Proposition \ref{c27Jan12}.(3))  and statement 2 imply the first two equalities in statement 1. Now, using Lemma \ref{b27Jan12},
\begin{eqnarray*}
 \TAut_K(P_n)_n\cap \CF_n &=& \Fix_{\TAut_K(P_n)_n} (\der_1, \ldots , \der_n) =\TAut_K(P_n)_n \cap \Fix_{\Aut_K(P_n)} (\der_1, \ldots , \der_n)\\
 &=& \TAut_K(P_n)_n\cap \Sh_n = \Sh_{n-1}.
 \end{eqnarray*}

 3. Since $\CT_n\times_{ex}\CF_n\subseteq \CU_n$ and $ \CU_n\cap \mT^n = \{ e\}$, statement 3 follows from statement 2:
 $$ \CU_n=\CU_n\cap G_n= \CU_n\cap (\mT^n \ltimes (\CT_n\times_{ex}\CF_n))= (\CU_n\cap \mT^n ) \ltimes (\CT_n\times_{ex}\CF_n)=\CT_n\times_{ex}\CF_n.\;\;\; \Box$$


\section{The group of automorphisms of the  Lie algebra $\ggu_n$}\label{GAUTUN}

The aim of this section is to find the groups $\CF_n$ (Theorem \ref{11Feb12}.(1)) and  $G_n$ (Theorem \ref{5Feb12}).

{\bf The $\ggu_n$-module $P_n$}. For each $n\geq 2$, $\ggu_n$ is a Lie subalgebra of the Lie algebra $\ggu_{n+1} = \ggu_n\oplus P_n\der_{n+1}$, $P_n\der_{n+1}$ is an ideal of the Lie algebra $\ggu_{n+1}$ and $[P_n\der_{n+1},P_n\der_{n+1}]=0$. In particular, $P_n\der_{n+1}$ is a left $\ggu_n$-module where the action of the Lie algebra $\ggu_n$ on $P_n\der_{n+1}$ is given by the rule (the adjoint action): $uv:= [u,v]$ for all $u\in \ggu_n$ and $v\in P_n\der_{n+1}$. The polynomial algebra $P_n$ is a left $\ggu_n$-module.

\begin{lemma}\label{aa21Jan12}
{\rm (Lemma 3.10, \cite{Lie-Un-GEN})}
\begin{enumerate}
\item The $K$-linear map $P_n\ra P_n\der_{n+1}$, $p\mapsto p \der_{n+1}$, is a $\ggu_n$-module isomorphism.
\item The $\ggu_n$-module $P_n$ is an indecomposable, uniserial $\ggu_n$-module, $\udim (P_n)=\o^n$ and $\ann_{\ggu_n}(P_n)=0$.
\end{enumerate}
\end{lemma}


The next proposition describes the algebra of all the $\ggu_n$-homomorphisms (and its group of units) of the $\ggu_n$-module $P_n$.  The $K$-derivation $\frac{d}{dx_n}$ of the polynomial algebra $K[x_n]$ is also denoted by $\der_n$.

\begin{proposition}\label{A26Jan12}
{\rm (Proposition 3.16, \cite{Lie-Un-GEN})}
\begin{enumerate}
\item The map $\End_{\ggu_n}(P_n) \ra \End_{K[\der_n]}(K[x_n]) = K[[\frac{d}{d x_n}]]$, $\v \mapsto \v |_{K[x_n]}$, is a $K$-algebra isomorphism with the inverse map $\v'\mapsto \v$ where $\v (x^\beta x_n^i)= X_{\beta , n} \v'(\frac{x_n^{i+1}}{i+1})$ for all $\beta \in \N^{n-1}$ and $ i\in \N$.
\item The map $\Aut_{\ggu_n} (P_n) \ra \Aut_{K[\der_n]}(K[x_n])= K[[\frac{d}{dx_n}]]^*$, $\v \mapsto \v |_{K[x_n]}$, is a group  isomorphism with the inverse map as in the statement 1  (where $K[[\frac{d}{dx_n}]]^*$ is the group of units of the algebra $K[[\frac{d}{dx_n}]]$).
\end{enumerate}
\end{proposition}

The partial derivatives $\der_1, \ldots , \der_n$ are commuting locally nilpotent derivations of the polynomial algebra  $P_n$.
 So, we can consider the skew power series algebra $P_n[[\der_1, \ldots , \der_n ; \frac{\der}{\der x_1},\ldots , \frac{\der}{\der x_n}]]$ which is also written as
 $P_n[[\der_1, \ldots , \der_n]]$, for short. Every element $s$ of this algebra is a unique (formal) series $\sum_{\alpha \in \N^n}p_\alpha \der^\alpha$ where $p_\alpha \in P_n$. The addition of two power series is defined in the obvious way $\sum_{\alpha \in \N^n}p_\alpha \der^\alpha +\sum_{\alpha \in \N^n}q_\alpha \der^\alpha =\sum_{\alpha \in \N^n}(p_\alpha +q_\alpha )  \der^\alpha$, and the multiplication satisfies the relations: $\der_ip=p\der_i+\der_i(p)$ for all $i=1, \ldots , n$ and $p\in P_n$. As the partial derivatives act locally nilpotently on the polynomial algebra $P_n$, the product of two power series can be written in the canonical form using the relations above, i.e.,
 $(\sum_{\alpha \in \N^n}p_\alpha \der^\alpha )(\sum_{\beta \in \N^n}q_\beta \der^\beta )=\sum_{\g \in \N^n}r_\g  \der^\g$ for some polynomials $r_\g\in P_n$.
 For every  series $s=\sum_{\alpha \in \N^n} p_\alpha \der^\alpha\in P_n [[\der_1, \ldots , \der_n]]$ (where $p_\alpha \in P_n$), the action $s*p=\sum_{\alpha \in \N^n}p_\alpha \der^\alpha * p$ is well-defined  (the infinite sum is in fact a finite one). The algebra $P_n [[\der_1, \ldots , \der_n]]$  is the completion of the Weyl algebra $A_n= P_n[\der_1, \ldots , \der_n]=\oplus_{\alpha \in \N^n}P_n\der^\alpha$. The polynomial algebra $P_n$ is a left $P_n [[\der_1, \ldots , \der_n]]$-module. It is easy to show that the algebra homomorphism $P_n [[\der_1, \ldots , \der_n]]\ra \End_K(P_n)$ is a monomorphism, and we identify the algebra $P_n [[\der_1, \ldots , \der_n]]$ with its image in $\End_K(P_n)$. The polynomial algebra $P_n=\cup_{i\geq 0} P_{n, \leq i}$ is a filtered algebra by the total degree of the variables where $P_{n, \leq i} := \sum\{ Kx^\alpha \, | \, |\alpha |:=\alpha_1+\cdots + \alpha_n \leq i\}$ ($P_{n, \leq i}P_{n, \leq j}\subseteq P_{n, \leq i+j}$ for all $i,j\geq 0$). The vector space  $$\End_{\deg} (P_n):= \{ f\in \End_K(P_n) \, | \, f(P_{n, \leq i})\subseteq P_{n, \leq i}\}$$ is a subalgebra of $\End_K(P_n)$. Let $\Aut_{K[\der_1, \ldots , \der_n]}(P_n)$ be the group of invertible $K[\der_1, \ldots , \der_n]$-endomorphisms of the $K[\der_1, \ldots , \der_n]$-module $P_n$ and $K[[\der_1, \ldots , \der_n]]^*$ be the group of units of the algebra $K[[\der_1, \ldots , \der_n]]$.

\begin{proposition}\label{A3Feb12}

\begin{enumerate}
\item $\End_K(P_n) = P_n[[\der_1, \ldots , \der_n]]$.
\item $\End_{\deg}(P_n)=P_n[[\der_1, \ldots , \der_n]]_{\deg}:= \{ \sum_{\alpha \in \N^n} p_\alpha \der^\alpha \, | \, \deg (p_\alpha ) \leq | \alpha |\} $.
\item  $\End_{K[\der_1, \ldots , \der_n]}(P_n) = K[[\der_1, \ldots , \der_n]]$.
    \item  $\Aut_{K[\der_1, \ldots , \der_n]}(P_n) = K[[\der_1, \ldots , \der_n]]^*$.
\end{enumerate}
\end{proposition}

{\it Proof}. 1. This is   well-known and easy to prove.

2. Let $R$ be the RHS of the equality in statement 2. The inclusion $\End_{\deg}(P_n) \supseteq R$ is obvious. We have to show that the reverse inclusion  holds. Let $s=\sum p_\alpha \der^\alpha \in \End_{\deg}(P_n)$. We have to prove that $\deg (p_\alpha ) \leq |\alpha |$ for all $\alpha$. We use induction on $d=|\alpha |$. The initial case $d=0$ is obvious as $p_0= s*1\in K$. Let $d>0$ and suppose that the statement holds for all $d'<d$. Fix $\alpha \in \N^n$ such that $|\alpha | = d$. Then
$$ P_{n,\leq d} \ni s*x^\alpha = \alpha ! p_\alpha +\sum_{|\beta | <d} p_\beta \der^\beta *x^\alpha ,$$
and so $p_\alpha \in P_{n, \leq d}$ (since $\sum_{|\beta | <d} p_\beta \der^\beta *x^\alpha\in P_{n , \leq d}$), as required.

3. An element $s=\sum_{\alpha \in \N^n} p_\alpha \der^\alpha\in \End_K(P_n)$ belongs to $\End_{K[\der_1, \ldots , \der_n]}(P_n)$ iff $[\der_i, s]=0$ for all $i=1, \ldots , n$ iff every  $p_\alpha\in \cap_{i=1}^n \ker_{P_n} (\der_i) = K$ for all $i=1, \ldots , n$  (since $ [\der_i , s] = \sum \frac{\der p_\alpha }{\der x_i}\der^\alpha )$ iff $s\in K[[\der_1, \ldots , \der_n]]$.

4. Statement 4 follows from statement 3. $\Box $


By Theorem 2.2, \cite{Bav-inform}, the $K$-algebra endomorphisms  $\phi_i:= \sum_{k\geq 0} (-1)^k \frac{x_i^k}{k!}\der_i^k:P_n\ra P_n$ (where $i=1, \ldots , n$) commute and their composition
\begin{equation}\label{PHIdef}
\phi := \prod_{i=1}^n \phi_i : P_n\ra P_n, \;\; \sum_{\alpha \in \N^n}\l_\alpha x^\alpha \mapsto \l_0,
\end{equation}
is the projection onto the field $K$ in the decomposition $P_n= \oplus_{\alpha \in \N^n} Kx^\alpha$. The next proposition is an easy corollary of this fact. \begin{proposition}\label{b14Feb12}
Let $s=\sum_{\alpha \in \N^n} \l_\alpha \der^\alpha \in K[[\der_1, \ldots , \der_n]]=\End_{K[\der_1, \ldots , \der_n]}(P_n)$ where $\l_\alpha \in K$. Then $\l_\alpha = \phi s(\frac{x^\alpha}{\alpha !})$ for all $\alpha \in \N^n$ where  $\alpha ! := \prod_{i=1}^n \alpha_i!$.
\end{proposition}

{\it Proof}. $\phi s(\frac{x^\alpha}{\alpha !})=\phi (\l_\alpha +\cdots ) = \l_\alpha$ where the three dots denote an element of the vector space $\oplus_{\alpha \in \N^n\backslash \{ 0\}}Kx^\alpha$. $\Box $


\begin{lemma}\label{a9Feb12}
Let $\s \in \CF_n$.
\begin{enumerate}
\item Then $\s |_{P_{n-1}\der_n} : P_{n-1}\der_n\ra P_{n-1}\der_n$ belongs to $I = 1+\sum_{i=1}^{n-1}\d_i K[[\d_1, \ldots , \d_{n-1}]]$ where $\d_i := \ad (\der_i)$ for  $i=1, \ldots , n-1$.
\item For all $\alpha \in \N^{n-1}\backslash \{ 0\}$, $\s (x^\alpha \der_n) = p_\alpha \der_n$ for some polynomial $p_\alpha \in P_{n-1}$ such that $p_\alpha = x^\alpha + \sum_{\beta \in \N^{n-1}, \beta\prec\alpha}\l_\beta x^\beta$ where $\beta = (\beta_i) \prec\alpha = (\alpha_i) $ iff $\beta_i\leq \alpha_i$ for all $i=1, \ldots , n-1$ and $\beta_j<\alpha_j$ for some $j$.
\end{enumerate}
\end{lemma}

{\it Proof}. 1. Recall that all ideals of the Lie algebra $\ggu_n$ are characteristic ideals, and the vector space $P_{n-1}\der_n$ is an ideal of the Lie algebra $\ggu_n$. Therefore, the restriction map $\tau  := \s |_{P_{n-1}\der_n}$ is a well-defined map.  Since $\s \in \CF_n$, the map $\tau$ commute with the inner derivations $\d_i $ where $i=1,\ldots , n-1$. The $\ggu_{n-1}$-modules  $P_{n-1}$ and $P_{n-1}\der_n$ are isomorphic (Corollary \ref{aa21Jan12}.(1)), hence $\tau \in \Aut_{K[\d_1, \ldots , \d_{n-1}]}(P_{n-1}\der_n)\simeq \Aut_{K[\der_1, \ldots , \der_{n-1}]}(P_{n-1})\simeq K[[\der_1, \ldots , \der_{n-1}]]^*$. Then,
$$\tau \in \Aut_{K[\d_1, \ldots , \d_{n-1}]}(P_{n-1}\der_n) \simeq K[[\d_1, \ldots , \d_{n-1}]]^*.$$ Since $\tau (\der_n)=\s (\der_n) = \der_n$,   we must have $\tau \in I$.

2. Statement 2 follows from statement 1. $\Box $


{\bf The subgroup  $\mF_n$ of $\CF_n$}.
Let $\ggu_1:= K\der_1$, the abelian 1-dimensional Lie algebra. For $n\geq 2$, the set
\begin{equation}\label{bFndef}
\mF_n:= \Fix_{\CF_n}(\ggu_{n-1}) = \{ \s \in \CF_n \, | \, \s (u) = u\;\; {\rm for \; all}\;\; u\in \ggu_{n-1}\}
\end{equation}
is a subgroup of $\CF_n$. Notice that $\mF_2= \CF_2$. Recall that every ideal of the Lie algebra $\ggu_n$ is a characteristic ideal (Corollary \ref{c10Dec11}), $P_{n-1}\der_n$ is an ideal of the Lie algebra $\ggu_n = \ggu_{n-1}\oplus P_{n-1}\der_n$ such that $\ggu_n / P_{n-1}\der_n \simeq \ggu_{n-1}$. In view of this Lie algebra isomorphism, for each natural number $n\geq 3$,  there is the group homomorphism
\begin{equation}\label{chindef}
\chi_n : \CF_n \ra \CF_{n-1}, \;\; \s \mapsto (u+P_{n-1}\der_n\mapsto \s (u) +P_{n-1}\der_n).
\end{equation}
It is obvious that $\mF_n\subseteq \ker (\chi_n)$. The ideal $P_{n-1}\der_n$ of the Lie algebra $\ggu_n$ is a (left) $\ggu_n$-module and a $\ggu_{n-1}$-module since $\ggu_{n-1}\subseteq \ggu_n$. Let $\Aut_{\ggu_{n-1}}(P_{n-1}\der_n)$ be the group of automorphisms of the $\ggu_{n-1}$-module $P_{n-1}\der_n$. The set $1+\der_{n-1}K[[\der_{n-1}]]$ is a subgroup of the group $K[[\der_{n-1}]]^*$ of units of the power series algebra $K[[\der_{n-1}]]$.

The following proposition is an explicit description of the group $\mF_n$.


\begin{proposition}\label{c28Jan12}
For $n\geq 2$, the map $\eta_n : 1+\der_{n-1} K[[\der_{n-1}]]\ra \mF_n$, $s= \sum \l_i \der_{n-1}^i\mapsto \eta_n(s)$, is a group isomorphism  where,  for $p\in P_{n-1}$,   $ \eta_n(s) :p\der_n\mapsto (\sum \l_i \frac{\der^i p }{\der x_{n-1}^i}) \der_n$.   So, $\eta_n(s)$ acts on the elements of the algebra  $\ggu_n$ as $\sum \l_i \ad ( \der_{n-1})^i$. In particular, the group $\mF_n=1+\d_{n-1}K[[\d_{n-1}]]$ is abelian where $\d_{n-1}= \ad (\der_{n-1})$ (equally, we can write $\mF_n = 1+\der_{n-1}K[[\der_{n-1}]]$). Moreover, $\eta_n = (\xi_{n-1}\pi_{n-1}\rho_n)^{-1}$ where the maps $\rho_n$, $\pi_{n-1}$ and $\xi_{n-1}$ are defined in the proof (see (\ref{FnA}), (\ref{FnaA1}) and (\ref{FnA2}) respectively).
\end{proposition}
{\it Proof}. Since $\mF_n = \Fix_{\CF_n} (\ggu_{n-1})$ and $P_{n-1}\der_n$ is an abelian characteristic ideal of the Lie algebra $\ggu_n$, the restriction map
\begin{equation}\label{FnA}
\rho_n: \mF_n\ra \Aut_{\ggu_{n-1}}(P_{n-1}\der_n)_n , \;\; \s \mapsto \s |_{P_{n-1}\der_n},
\end{equation}
is a group isomorphism where $\Aut_{\ggu_{n-1}}(P_{n-1}\der_n)_n := \{ \v \in \Aut_{\ggu_{n-1}}(P_{n-1}\der_n)\, | \, \v (\der_n ) = \der_n \}$. The map
\begin{equation}\label{FnA1}
P_{n-1}\der_n \ra P_{n-1} , \;\; p\der_n \mapsto p,
\end{equation}
is a $\ggu_{n-1}$-module isomorphism, and so $\Aut_{\ggu_{n-1}}(P_{n-1}\der_n) \simeq \Aut_{\ggu_{n-1}}(P_{n-1})$ and
\begin{equation}\label{FnaA1}
\pi_{n-1}: \Aut_{\ggu_{n-1}}(P_{n-1}\der_n)_n\simeq  \Aut_{\ggu_{n-1}}(P_{n-1})_1:= \{ \s \in \Aut_{\ggu_{n-1}}(P_{n-1})\, | \, \s (1) = 1 \}.
\end{equation}
By Proposition \ref{A26Jan12}.(1,2), the  restriction map
\begin{equation}\label{FnA2}
\xi_{n-1}: \Aut_{\ggu_{n-1}}(P_{n-1})_1\ra \Aut_{K[\der_{n-1}]} (K[x_{n-1}])_1= 1+\der_{n-1}K[[\der_{n-1}]],
  \;\; \v \mapsto \v |_{K[x_{n-1}]},
\end{equation}
is a group isomorphism where $\Aut_{K[\der_{n-1}]}(K[x_{n-1}])_1: = \{ \s \in \Aut_{K[\der_{n-1}]} (K[x_{n-1}])\, | \, \s (1) = 1\}$. Then the automorphism $\eta_n$ is equal to $(\xi_{n-1}\pi_{n-1}\rho_n)^{-1}$. $\Box $


{\it Remark}. It is useful to identify the groups $\mF_n = 1+\d_{n-1}K[[\d_{n-1}]]$  and $1+\der_{n-1}K[[\der_{n-1}]]$ via the isomorphism $\eta_n$, i.e.,
\begin{equation}\label{Fn=idnt}
\mF_n = 1+\d_{n-1}K[[\d_{n-1}]]=1+\der_{n-1}K[[\der_{n-1}]].
\end{equation}
The first equality is used when the action of automorphisms of the group $\mF_n$ on elements of the Lie algebra $\ggu_n$ is considered. The second equality is used when we want to stress how the polynomial coefficient of $\der_n$   is  changed under this action on the ideal $P_{n-1}\der_n$ of the Lie algebra $\ggu_n$.


{\bf The subgroup $\mE_n$ of $\CF_n$}.
For $n\geq 3$, let
\begin{eqnarray}\label{En=Fix}
 \mE_n &:=& \Fix_{\ker (\chi_n)}(P_{n-1}\der_n) = \{ \s \in \ker (\chi_n) \, | \, \s (u) = u\;\; {\rm for \; all}\;\; u\in P_{n-1}\der_n\}\\
 &=&\{ \s \in \CF_n \, | \, \s (u)=u \;\; {\rm for \; all}\;\; u\in P_{n-1}\der_n; \s (v) -v\in P_{n-1}\der_n\; {\rm for\; all}\; v\in \ggu_{n-1}\}.
\end{eqnarray}
Clearly, $\mF_n \cap \mE_n = \{ e\}$ since $\ggu_n = \ggu_{n-1}\oplus P_{n-1} \der_n$. Consider the vector space (of certain 1-cocycles)
 $ Z^1_{n-1} := \{ c\in \Hom_K(\ggu_{n-1} , P_{n-1}\der_n)\, | \, c(\der_1) = \cdots = c(\der_{n-1})=0$ and
\begin{equation}\label{cuv}
c([u,v]) = [c(u), v]+[u,c(v)]
\end{equation}
for all $u,v\in \ggu_{n-1}\}$. In particular, $Z^1_{n-1}$ is an additive (abelian) group.

\begin{lemma}\label{a4Feb12}
For $n\geq 3$, the map
$\D_n : \mE_n\ra Z^1_{n-1}, \s \mapsto \s -1 : u\mapsto \s (u) - u$
(where $u\in \ggu_{n-1}$),   is a group isomorphism with the inverse map $c\mapsto \s_c$ where $\s_c (u) = \begin{cases}
u+c(u)& \text{if }u\in \ggu_{n-1},\\
u& \text{if }u\in P_{n-1}\der_n.\\
\end{cases}$ In particular, $\mE_n$ is an abelian group  (and a vector space over the field $K$).
\end{lemma}

{\it Proof}. Let us show that the map $\D_n$ is well-defined. Let $\s \in \mE_n$. Then $ c = \s -1\in \Hom_K(\ggu_{n-1} , P_{n-1}\der_n )$ (since $ \s \in \ker (\chi_n)$), and $c(\der_1) = \cdots = c(\der_n) =0$ (since  $\s \in \CF_n$). The condition (\ref{cuv}) follows from the facts that $\s $ is a Lie algebra homomorphism, $\im ( c) \subseteq P_{n-1}\der_n$ and $[P_{n-1}\der_n , P_{n-1}\der_n]=0$. In more details, for all elements $u,v\in \ggu_{n-1}$,
\begin{eqnarray*}
0&=& \s ([u,v]) -[\s (u) , \s (v)]=[u,v]+c([u,v]) - [u,v]-[c(u), v]- [u,c(v)]-[c(u), c(v)] \\
 &=& c([u,v])-[c(u), v]-[u,c(v)].
\end{eqnarray*}
 The map $\D_n$ is a group homomorphism: for all elements $\s_1, \s_2\in \mE_n$,
$$
\D_n(\s_1\s_2) =\s_1\s_2-1= \s_1-1+\s_1(\s_2-1) = \s_1-1+\s_2-1= \D_n(\s_1) + \D_n( \s_2),$$
we used the fact that $\im (\s_2 -1) \subseteq P_{n-1}\der_n$ and $\s_1(u) = u$ for all elements $u\in P_{n-1}\der_n$ (since $\s_1\in \mE_n$).  Since $\ker( \D_n) = \{ \id \}$, the map $\D_n$ is a monomorphism.

To finish the proof of the lemma it suffices to show that the map $Z^1_{n-1}\ra \mE_n$, $ c\mapsto \s_c$, is well-defined (indeed, this claim guarantees the surjectivity of the map $\D_n$, i.e.,   $\D_n$  is a group isomorphism; then it is obvious that the map $c\mapsto \s_c$ is the inverse of $\D_n$). By the very definition, the $K$-linear map $\s_c:\ggu_n\ra \ggu_n$ is a bijection. To finish the proof of the claim it suffices to show that $\s_c$ is a Lie algebra homomorphism (since then $\s_c\in \ker (\chi_n)$ and $\s_c\in \Fix_{\ker (\chi_n)} ( P_{n-1}\der_n)= \mE_n$). Since $\ggu_n = \ggu_{n-1}\oplus P_{n-1}\der_n$ and $[P_{n-1}\der_n , P_{n-1}\der_n ]=0$, it suffices to verify that, for all elements $u,v\in \ggu_{n-1}$, $\s_c([u,v]) = [ \s_c(u),  \s_c(v)]$. This follows from (\ref{cuv}) and $[c(u), c(v)]=0$:
$$\s_c([u,v])-  [ \s_c(u),  \s_c(v)]=[u,v]+c([u,v])-[u,v]-
 [ c(u), v]-[u, c(v)]-[c(u), c(v)]=0.\; \Box$$


{\bf The subgroup $\ker (\chi_n )$ of $\CF_n$}.
The next corollary describes the groups $\ker (\chi_n)$ when $n\geq 3$.
\begin{corollary}\label{b4Feb12}
Let $n\geq 3$. Then
\begin{enumerate}
\item $\ker(\chi_n) = \mF_n\ltimes \mE_n$.
\item Each element $\s \in \ker (\chi_n)$ is the unique product $ef$ where $e\in \mE_n$,  $e(u) = u+c(u)$ and $c(u):= \s (u) - u$  for all $u\in \ggu_{n-1}$, and $f:=e^{-1}\s\in \mF_n$.
\end{enumerate}
\end{corollary}

{\it Proof}. It is obvious that $\ \mF_n \cap \mE_n = \{ e\}$ and $\tau \mE_n \tau^{-1} \subseteq \mE_n$ for all $\tau \in \mF_n$. Therefore, $\mF_n \ltimes \mE_n\subseteq \ker (\chi_n)$. Now, to finish the proof of both statements it suffices to show that each element $\s \in \ker (\chi_n)$ is the product $ef$ where $e$ and $f$ are as in statement 2.  It is easy to check that $c\in Z^1_{n-1}$ (where $ c(u) = \s (u) - u$ for $u\in \ggu_{n-1}$): for all elements $u,v\in \ggu_{n-1}$,
\begin{eqnarray*}
 c([u,v])&=& \s ([u,v]) -[u,v]= [\s (u), \s (v)]-[u,v]= [u+c(u), v+c(v)]-[u,v]\\
 &=& [c(u), v]+[u,c(v)].
\end{eqnarray*}
By Lemma \ref{a4Feb12}, $\D_n^{-1}(c)  \in \mE_n$. Notice that $\D_n^{-1} (c) = e$. Now, $e^{-1}\s \in \mF_n$ since, for all elements $u\in \ggu_{n-1}$,
$$ e^{-1}\s (u) = e^{-1} (u+c(u)) = e^{-1}(u) +c(u) = u-c(u) + c(u) = u.\;\;\;\; \Box$$


We will see that $\ker (\chi_n) = \mF_n\times \mE_n$ (Theorem \ref{11Feb12}.(3)).

Let $n\geq 3$. For each $i=1,\ldots , n-2$, $\ggu_i+\CD_{n-1}$ is a Lie subalgebra of the Lie algebra $\ggu_n$ where $\ggu_1:=K\der_1$ and $\CD_{n-1}:= \oplus_{i=1}^{n-1}K\der_i$. Notice that the Lie algebra $\ggu_n$ is the (adjoint) $\ggu_i+\CD_{n-1}$-module, and the vector spaces $P_i\der_{i+1}$, $P_i\der_n$ and $P_{n-1}\der_n$ are $\ggu_i+\CD_{n-1}$-submodules of $\ggu_n$.  The $\ggu_i+\CD_{n-1}$-modules $P_i\der_{i+1}$ and $P_i\der_n$ are annihilated by the elements $\der_{i+1}, \ldots , \der_{n-1}$.
 Our goal is to give  an explicit description of the group $Z_{n-1}^1$ (Corollary \ref{d4Feb12}). The group $Z_{n-1}^1$ turns out to be the direct product of certain subgroups described in Proposition \ref{c4Feb12}.
 Let
 \begin{eqnarray*}
 \Hom_{\ggu_i+\CD_{n-1}}(P_i\der_{i+1}, P_{n-1}\der_n)_0&: =&\{ \v \in \Hom_{\ggu_i+\CD_{n-1}}(P_i\der_{i+1}, P_{n-1}\der_n)\, | \, \v (\der_{i+1})=0\}, \\
 \Hom_{\ggu_i+\CD_{n-1}}(P_i\der_{i+1}, P_i\der_n)_0&: =&\{ \v \in \Hom_{\ggu_i+\CD_{n-1}}(P_i\der_{i+1}, P_i\der_n)\, | \, \v (\der_{i+1})=0\}.
\end{eqnarray*}
\begin{proposition}\label{c4Feb12}
Let $n\geq 3$ and $\ggu_1= K\der_1$. Then
\begin{enumerate}
\item For all $i=1, \ldots , n-2$, $\Hom_{\ggu_i+\CD_{n-1}}(P_i\der_{i+1}, P_{n-1}\der_n) = \Hom_{\ggu_i+\CD_{n-1}}(P_i\der_{i+1}, P_i\der_n) \simeq \End_{\ggu_i}(P_i)\simeq K[[t]]$. Moreover, the map
    $$\alpha_n : K[[t]]\ra \Hom_{\ggu_i+\CD_{n-1}}(P_i\der_{i+1}, P_{n-1}\der_n), \;\; \sum_{j\geq 0} \l_j t^j \mapsto \v, $$
(where $\l_j\in K$) is an isomorphism of vector spaces where, for all elements $\beta \in \N^{i-1}$ and $k\in \N$,
    $$ \v ( x^\beta x_i^k\der_{i+1}):= [ X_{\beta , i }, \sum_{j\geq 0} \l_j (\ad \, \der_i)^j ((k+1)^{-1} x_i^{k+1}\der_n)]=\sum_{j\geq 0}\l_j(\ad \, \der_i)^j(x^\beta x_i^k\der_n).$$
\item For all $i=1, \ldots , n-2$, $\Hom_{\ggu_i+\CD_{n-1}}(P_i\der_{i+1}, P_{n-1}\der_n)_0 = \Hom_{\ggu_i+\CD_{n-1}}(P_i\der_{i+1}, P_i\der_n)_0 \simeq tK[[t]]$. Moreover, the map
    $$\alpha_n : tK[[t]]\ra \Hom_{\ggu_i+\CD_{n-1}}(P_i\der_{i+1}, P_{n-1}\der_n)_0, \;\; \sum_{j\geq 1} \l_j t^j \mapsto \v, $$
(where $\l_j\in K$) is an isomorphism of vector spaces where, for all elements $\beta \in \N^{i-1}$ and $k\in \N$,
    $$ \v ( x^\beta x_i^k\der_{i+1}):= [ X_{\beta , i }, \sum_{j\geq 1} \l_j (\ad \, \der_i)^j ((k+1)^{-1} x_i^{k+1}\der_n)]=\sum_{j\geq 1}\l_j(\ad \, \der_i)^j(x^\beta x_i^k\der_n).$$
    \item For every $i=1, \ldots , n-2$, the $K$-linear map
$$ \beta_{n,i} : \Hom_{\ggu_i+\CD_{n-1}}(P_i\der_{i+1}, P_{n-1}\der_n)_0\ra Z^1_{n-1}, \;\; \psi \mapsto c_\psi,$$ is an injection where, for $u\in P_{j-1}\der_j$, $j=1, \ldots , n-1$, $c_\psi (u) := \begin{cases}
\psi (u) & \text{if }j=i+1,\\
0& \text{if }j\neq i+1.\\
\end{cases}$
\end{enumerate}
\end{proposition}

{\it Proof}. 1. Let $\d_i := \ad (\der_i)$ for $i=1, \ldots , n-1$ and  $\psi \in H:= \Hom_{\ggu_i+\CD_{n-1}}(P_i\der_{i+1}, P_{n-1}\der_n)$.
 For all $j=i+1, \ldots , n-1$, $[\der_j , \psi (P_i\der_{i+1})]= \psi ([\der_j, P_i\der_{i+1}]) = \psi (0) =0$, hence $$\im (\psi ) \subseteq \bigcap_{j=i+1}^{n-1}\ker_{P_{n-1}\der_n} (\d_j)=(\bigcap_{j=i+1}^{n-1}\ker_{P_{n-1}} (\der_j))\der_n=P_i\der_n.$$ Therefore, $H= \Hom_{\ggu_i+\CD_{n-1}}(P_i\der_{i+1}, P_i\der_n)$. The maps $$P_i\der_{i+1}\ra P_i\der_n\ra P_i,\;\;  p\der_{i+1}\ra p\der_n\ra p$$ (where $p\in P_i$) are $\ggu_i+\CD_{n-1}$-module isomorphisms. The $\ggu_i+\CD_{n-1}$-modules $P_i\der_{i+1}$ and $P_i\der_n$ are annihilated by the elements $\der_{i+1}, \ldots , \der_{n-1}$. So,
 $$H= \Hom_{\ggu_i+\CD_{n-1}}(P_i\der_{i+1}, P_i\der_n)\simeq \End_{\ggu_i+\CD_{n-1}}(P_i) \simeq \End_{\ggu_i}(P_i)\simeq K[[t]],$$
by Proposition \ref{A26Jan12}.(1). The map $\alpha_n: K[[t]]\ra H$ is the inverse of the above isomorphism $H\simeq K[[t]]$  (see Proposition \ref{A26Jan12}.(1)).

2. Statement 2 follows from statement 1.

3. Let $c=c_\psi$. By the very definition of $c$, $c(\der_1) = \cdots = c(\der_{n-1})=0$. We have to show that $c([u,v]) = [c(u), v]+[u,c(v)]$ for all elements $u\in P_{s-1}\der_s$ and $v\in P_{t-1}\der_t$ where $s,t=1, \ldots , n-1$. Without loss of generality we may assume that $s\leq t$.

If $s\neq i+1$ and $t\neq i+1$ then $[u,v]\not\in P_i\der_{i+1}$, and the equality above trivially holds ($0=0+0$).

If  $s< i+1$ and $t= i+1$ then $c(u)=0$, and the equality that we have to check reduces to the equality $ \psi ( [u,v]) = [u , \psi (v)]$ which is obviously true as the map $\psi$ is a $\ggu_i+\CD_{n-1}$-homomorphism.

If $s=t=i+1$ then the equality that we have to check reduces to the equality $\psi ([u,v])=[\psi (u), v]+[u,\psi (v)]$ which is obviously true as $[u,v]\in [P_i\der_{i+1}, P_i\der_{i+1}]=0$, $[\psi (u), v] \in [P_i\der_n, P_i\der_{i+1}]=0$, and $[u,\psi (v)]\in [ P_i\der_{i+1}, P_i\der_n]=0$ (since $\im (\psi ) \subseteq P_i\der_n$, see the proof of statement 1).  $\Box $


In combination  with Proposition \ref{c4Feb12}, the following
 corollary gives an explicit description of the vector space $Z^1_{n-1}$, and, as a result, using Lemma \ref{a4Feb12} we have an  explicit description
 of the group $\mE_n$ and its generators.

\begin{corollary}\label{d4Feb12}
Let $n\geq 3$. Then the $K$-linear map
$$ \beta_n:=\bigoplus_{i=1}^{n-2}\beta_{n,i}:\bigoplus_{i=1}^{n-2}
\Hom_{\ggu_i+\CD_{n-1}}(P_i\der_{i+1}, P_{n-1}\der_n)_0\ra Z^1_{n-1}, \;\; (\psi_1, \ldots , \psi_{n-2})\mapsto c_{\psi_1}+\cdots +c_{\psi_{n-2}},$$
 is a bijection where $c_{\psi_i}= \beta_{n,i}(\psi_i)$. In particular, $Z^1_{n-1}\simeq (tK[[t]])^{n-2}$, the direct sum of $n-2$ copies of the vector space $tK[[t]]$.
\end{corollary}

{\it Proof}. In view of the explicit nature of the maps $\beta_{n,i}$, the map $\beta_n$ is an injection. It remains to show that the map $\beta_n$ is a surjection. Let $c\in Z^1_{n-1}$. Then, by (\ref{cuv}),  $$\psi_1:= c|_{P_1\der_2}
\in \Hom_{\CD_{n-1}}(P_1\der_2, P_{n-1}\der_n)_0 = \Hom_{\ggu_1+\CD_{n-1}}(P_1\der_2, P_{n-1}\der_n)_0,$$ and so $c_{\psi_1} := \beta_{n,1}(\psi_1)\in Z^1_{n-1}$, by Proposition \ref{c4Feb12}.(3). By (\ref{cuv}), $$c_2':= c-c_{\psi_1}\in \Hom_{\ggu_2+\CD_{n-1}}(\ggu_{n-1}, P_{n-1}\der_n)_0$$ since $c_2'(\ggu_2+\CD_{n-1})=0$. Then $\psi_2:= c_2'|_{P_2\der_3}\in \Hom_{\ggu_2+\CD_{n-1}}(P_2\der_3, P_{n-1}\der_n)_0$, and so $c_{\psi_2}:= \beta_{n,2}(\psi_2) \in Z^1_{n-1}$,  by Proposition \ref{c4Feb12}.(3). Then $c_3':= c-c_{\psi_1}-c_{\psi_2}\in \Hom_{\ggu_3+\CD_{n-1}}(\ggu_{n-1}, P_{n-1}\der_n)_0$, by (\ref{cuv}) and the fact that $c_3'(\ggu_3+\CD_{n-1}) =0$. Continue this process (or use induction) we obtain the decomposition $c=c_{\psi_1}+\cdots + c_{\psi_{n-2}}$ where $c_{\psi_i}\in \im (\beta_{n,i})$ for all $i=1,\ldots , n-2$. This mean that the map $\beta_n$ is a surjective map, as required. $\Box $


{\bf The structure of the group $\mE_n$}.
For each $i=2, \ldots , n-1$, let $\mE_{n,i}:= \im (\D_n^{-1}\beta_{n, i-1})$ (notice the shift by 1 of the indices  when comparing them with the indices in Corollary \ref{d4Feb12}). Then $\mE_{n,i}= \{ e_i'(s_i)\, | \, s_i\in \der_{i-1} K[[\der_{i-1}]]\} \simeq (\der_{i-1} K[[\der_{i-1}]], +)$ (via $ e_i'(s_i) \mapsto s_i$)  where, for all $j=1, \ldots , n$ and $\alpha \in \N^{j-1}$,
\begin{equation}\label{E2n1}
 e_i'(s_i)(x^\alpha \der_j) = \begin{cases}
x^\alpha \der_j +s_i(x^\alpha ) \der_n& \text{if }j=i,\\
x^\alpha \der_j& \text{if }j\neq i, \\
\end{cases}
\end{equation}
\begin{equation}\label{E2n}
\mE_n = \prod_{i=2}^{n-1}\mE_{n,i}.
\end{equation}
So, each element $e'\in \mE_n$ is the unique product $e'=e'_2\cdots e_{n-1}'$ with $e_i' =e_i'(s_i)\in \mE_{n,i}$ and each automorphism $e_i'$ is uniquely determined by a series $s_i = \sum_{j\geq 1} \l_{ij} \der_{i-1}^j$  with $\l_{ij}\in K$. By (\ref{E2n1}), for all elements $x^\alpha \der_j \in \ggu_n$ where $\alpha \in \N^{j-1}$,
\begin{equation}\label{E2n2}
e'(x^\alpha \der_j)=\begin{cases}
e_j'(x^\alpha \der_j)& \text{if }2\leq j\leq n-1,\\
x^\alpha \der_j& \text{if }j=1,n, \\
\end{cases}
=\begin{cases}
x^\alpha \der_j+s_j(x^\alpha) \der_n& \text{if }2\leq j\leq n-1,\\
x^\alpha \der_j& \text{if }j=1,n. \\
\end{cases}
\end{equation}
Equivalently, for all elements $ u = \sum_{i=1}^n p_i\der_i \in \ggu_n$ where $p_i\in P_{i-1}$ for all $i=1, \ldots , n$,
\begin{equation}\label{E2n2a}
e'(u) = u+\sum_{i=2}^{n-1}s_i(p_i) \der_n.
\end{equation}

\begin{lemma}\label{c9Feb12}
For $n\geq 4$, $\im (\chi_n) \cap \mE_{n-1}=\{ e\}$.
\end{lemma}

{\it Proof}. Let $\s\in \mE_{n-1}$. Then for all polynomials  $p_{i}\in P_{i}$, $i=1, \ldots , n-2$, by (\ref{E2n2a}),
$$ \s ( \sum_{i=1}^{n-2} p_i\der_i) = \sum_{i=1}^{n-2} p_i\der_i +(\sum_{i=2}^{n-2}s_i(p_i))\der_{n-1}$$
where $ s_i =\sum_{j\geq 1} \l_{ij} (\ad \, \der_{i-1})^j$ for some scalars $\l_{ij}\in K$. Suppose that the automorphism $\s$ belongs to $\im (\chi_n)$, i.e.,  $\s = \chi_n( \s')$ for some $\s'\in \CF_n$. By Lemma \ref{a9Feb12}, $\s' ( x_{n-1}\der_n)= x_{n-1}\der_n +\l \der_n$ for some $\l \in K$. Then applying the automorphism $\s'$ to the equality $[\sum_{i=1}^{n-2}p_{i}\der_i, x_{n-1}\der_n] =0$ in the algebra $\ggu_n$  yields the equality $(\sum_{i=2}^{n-2}s_i (p_i) ) \der_n=0$ for all polynomials $p_i\in P_{i-1}$, $i=2, \ldots , n-2$. We used the fact that $\der_n$ is a central element of the Lie algebra $\ggu_n$ and $[P_{n-1}\der_n , P_{n-1}\der_n]=0$. Hence, all $s_i=0$, i.e.,  $\s =e$. $\Box $


By Lemma \ref{a9Feb12}.(1), there is the short exact sequence of group homomorphisms
\begin{equation}\label{resFnd}
1\ra \ker (\res_n) \ra \CF_n\stackrel{\res_n}{\ra} 1+\d_{n-1}K[[\d_{n-1}]]\ra 1,
\end{equation}
where $\res_n(\s ) := \s|_{K[x_{n-1}]\der_n}: K[x_{n-1}]\der_n\ra K[x_{n-1}]\der_n$ with $\res_n(\mF_n) =1+\d_{n-1}K[[\d_{n-1}]]$ (Proposition \ref{c28Jan12}), and the natural inclusion $\mF_n=1+\d_{n-1}K[[\d_{n-1}]]\subseteq \CF_n$ is a splitting for the epimorphism $\res_n$. Therefore,
\begin{equation}\label{FFn=K}
\CF_n= \mF_n\ltimes \ker (\res_n).
\end{equation}
The group $\CF_n$ contains the groups $\Sh_{n-2}$, $\mF_n$ and $\mE_n$. Let us show that  {\em the  elements of the groups $\Sh_{n-2}$, $\mF_n$ and $\mE_n$ pairwise commute.} Let $u=\sum_{i=1}^np_i\der_i = u_{n-1}+u_n$ where $p_i\in P_{i-1}$ for all $i=1, \ldots , n$, $u_{n-1}=\sum_{i=1}^{n-1}p_i\der_i$, $u_n=p_n\der_n$, $e'\in \mE_n$, $f=1+\sum_{i\geq 1} \l_i \der_{n-1}^i\in \mF_n$ where $\l_i\in K$ and $s\in \Sh_{n-1}$. We assume that (\ref{E2n2a}) holds for the element $e'$. Then
\begin{eqnarray*}
 e'f(u) &=& e'(u_{n-1}+f(u_n))= u_{n-1}+\sum_{i=2}^{n-1} s_i (p_i) \der_n + f(u_n),  \\
 fe'(u)&=& f(u+\sum_{i=2}^{n-1} s_i (p_i) \der_n )=u_{n-1}+f(u_n)+\sum_{i=2}^{n-1} s_i (p_i) \der_n.
\end{eqnarray*}
The last equality holds since $f(\sum_{i=2}^{n-1} s_i (p_i) \der_n )=\sum_{i=2}^{n-1} s_i (p_i) \der_n$. This follows from the inclusions  and $s_i(p_i) \in P_{i-1}$ for all $i=2, \ldots , n-1$ ($\der_{n-1}(s_i(p_i))=0$ for all $i=2, \ldots , n-1$). Therefore, $e'f=fe'$. Every element $s\in \Sh_{n-2}$ can be uniquely written as $s= e^{\sum_{i=1}^{n-2}\l_i\der_i}$ where $\l_i\in K$. Then it is obvious that $sf=fs$. Finally, $se =es$ since, for all elements $u\in \ggu_n$,
$$ se(u) = s(u+\sum_{i=2}^{n-1} s_i (p_i) \der_n)= s(u) +\sum_{i=2}^{n-1}  ss_i(p_i) \der_n
 = s(u) +\sum_{i=2}^{n-1} s_i s(p_i) \der_n=es(u), $$ as $ss_i = s_is$ for all elements $i=1, \ldots , n-1$, and $s(\der_j) = \der_j$ for all elements $j=1, \ldots , n$.

By Corollary \ref{b4Feb12}.(1), $\ker (\chi_n) = \mF_n\times \mE_n$. For $n=2$, $\Sh_0:=\{ e\}$ and $\mE_2:= \{ e\}$. The subgroup of the group $\CF_n$  that  these three groups generate is an abelian group. It is easy to see that $\Sh_{n-2}\cap ( \mF_n\times \mE_n)=\{ e\}$. Hence, $\Sh_{n-2}\times \mF_n\times \mE_n\subseteq \CF_n$. The next theorem shows that, in fact, the equality holds.
\begin{theorem}\label{11Feb12}

\begin{enumerate}
\item $\CF_n = \Sh_{n-2}\times \ker (\chi_n) = \Sh_{n-2}\times \mF_n\times \mE_n$.
\item $\im (\chi_n)= \Sh_{n-2}$.
\item The short exact sequence of group homomorphisms $1\ra \ker (\chi_n) \ra \CF_n\stackrel{\chi_n}{\ra}\Sh_{n-2}\ra 1$ is a split short exact sequence and the natural inclusion $\Sh_{n-2}\subseteq \CF_n$ is a splitting of the epimorphism $\chi_n$, $\ker (\chi_n) = \mF_n\times \mE_n$.
    \item $\ker (\res_n) = \Sh_{n-2}\times \mE_n$ and $\CF_n = \mF_n \times \ker (\res_n)$.
\end{enumerate}
\end{theorem}

{\it Proof}. 1-3. We use induction on $n\geq 2$ to prove all three  statements simultaneously. The initial step $n=2$ is trivial as $\CF_2=\mF_2$, $\CF_1=\{ e\}$, $\Sh_0:=\{ e\}$ and $\mE_1=\{ e\}$. So, let $n>2$ and we assume that all three statements hold for all $n'<n$. By induction, $\CF_{n-1}= \Sh_{n-3}\times \mF_{n-1}\times \mE_{n-1}$, and  as a result $$\chi_n : \CF_n\ra \CF_{n-1} = \Sh_{n-3}\times \mF_{n-1}\times \mE_{n-1}.$$ To show that statements 1-3 hold it suffices only  to prove statement 2, that is $\im (\chi_n) = \Sh_{n-2}$ (since then statement 3 follows as $\ker (\chi_n) = \mF_n \times \mE_n$; statement 3  and Corollary \ref{b4Feb12}.(1) imply statement 1). The subgroup $\Sh_{n-3}$ of $\CF_n$ is mapped isomorphically/identically onto its image $\chi_n(\Sh_{n-3}) = \Sh_{n-3}$. Then $\im (\chi_n) = \Sh_{n-3}\times I$ where $I:= \im (\chi_n) \cap (\mF_{n-1}\times \mE_{n-1})$. Notice that $\Sh_{n-2} = \Sh_{n-3} \times \sh_{n-2}$, $\sh_{n-2}= e^{K\der_{n-2}} \subseteq \mF_{n-1} = 1+\der_{n-2}K[[\der_{n-2}]]$ (recall the identification (\ref{Fn=idnt})). To finish the proof it remains to show that $I =\sh_{n-2}$.

 Let $\s \in \CF_n$ be such that $\s':=\chi_n( \s ) \in \mF_{n-1}\times \mE_{n-1}$, $\s'=fe'$ for some automorphisms $f\in \mF_{n-1}$ and $e'\in \mE_{n-1}$. We have to show that $\s'\in \sh_{n-2}$. Since $\mF_n\subseteq \ker (\chi_n)$, in view of (\ref{FFn=K}),  without loss of generality we may assume that $\s \in \ker (\res_n)$, that is $\s ( x_{n-1}^i\der_n) = x_{n-1}^i\der_n$ for all $i\geq 0$. Recall that $\ggu_n = \ggu_{n-2}\oplus P_{n-2}\der_{n-1}\oplus P_{n-1}\der_n= \ggu_{n-1}\oplus P_{n-1}\der_n$. For all elements $u\in P_{n-2}\der_{n-1}$,
 $$ \s ( u) = fe'(u) +c(u) = f(u)+ c(u)$$ for some map $c\in \Hom_K(P_{n-2}\der_{n-1}, P_{n-1}\der_n)$ ($e'(u) = u$ for all elements $u\in P_{n-2}\der_{n-1}$). The automorphism $f=1+\sum_{i\geq 1} \l_i\der_{n-2}^i\in \mF_{n-1} = 1+\der_{n-2} K[[\der_{n-2}]]$ (where $\l_i\in K$) acts on the elements $x^\alpha \der_i\in \ggu_{n-1}$ as $f(x^\alpha \der_i) = f(x^\alpha ) \der_i$ where
 $$ f(x^\alpha ) =  (1+\sum_{i\geq 1} \l_i\der_{n-2}^i) (x^\alpha ).$$
When we apply the automorphism $\s $ to the following identities in the Lie algebra $\ggu_n$,
$$ [ x_{n-2}^i\der_{n-1}, (j+1)^{-1} x_{n-1}^{j+1}\der_n]= x_{n-2}^ix_{n-1}^j\der_n, \;\;\; i,j\in \N, $$ it yields the identities
\begin{eqnarray*}
 \s (x_{n-2}^ix_{n-1}^j\der_n)&=& [ \s (x_{n-2}^i\der_{n-1}), \s ((j+1)^{-1} x_{n-1}^{j+1}\der_n)]\\
 &=&[ f(x_{n-2}^i)\der_{n-1}+c(x_{n-2}^i\der_{n-1}), (j+1)^{-1} x_{n-1}^{j+1}\der_n]
 = f(x_{n-2}^i) x_{n-1}^j\der_n
\end{eqnarray*}
since $[c(x_{n-2}^i\der_{n-1}), (j+1)^{-1} x_{n-1}^{j+1}\der_n] \subseteq [ P_{n-1}\der_n , P_{n-1}\der_n]=0$.  Similarly, when we apply  the automorphism $\s $ to the following  identities in the Lie algebra $\ggu_n$,
$$[x_{n-2}^i\der_{n-1}, x_{n-2}^jx_{n-1}\der_n] = x_{n-2}^{i+j}\der_n, \;\;\;  i,j\in \N,$$ we deduce the identities
$$ f(x_{n-2}^i)f(x_{n-2}^j)\der_n=f(x_{n-2}^{i+j})\der_n, \;\;\; i,j\in \N .$$
In more detail,
\begin{eqnarray*}
f(x_{n-2}^{i+j})\der_n &= & \s (x_{n-2}^{i+j}\der_n)=\s ([ x_{n-2}^i\der_{n-1}, x_{n-2}^jx_{n-1}\der_n])
=[\s (x_{n-2}^i\der_{n-1}), \s (x_{n-2}^jx_{n-1}\der_n)]  \\
 &=& [f(x_{n-2}^i\der_{n-1})+c(x_{n-2}^i\der_{n-1}), f(x_{n-2}^j) x_{n-1}\der_n]\\
 &=& [f(x_{n-2}^i)\der_{n-1}, f(x_{n-2}^j)x_{n-1}\der_n]=
  f(x_{n-2}^i)f(x_{n-2}^j)\der_n
\end{eqnarray*}
since $[c(x_{n-2}^i\der_{n-1}), f(x_{n-2}^j)x_{n-1}\der_n]\in [P_{n-1}\der_n, P_{n-1}\der_n]=0$. Therefore, $f(x_{n-2}^{i+j}) = f(x_{n-2}^i) f(x_{n-2}^j)$ for all $i,j\in \N$. This means that the map $f=1+\sum_{i\geq 1} \l_i \der_{n-2}^i: K[x_{n-2}]\ra K[x_{n-2}]$ is an automorphism of the polynomial algebra $K[x_{n-2}]$. Since $f(x_{n-2})= x_{n-2}+\l_1$, we must have $f= e^{\l_1\der_{n-2}}\in \sh_{n-2}$. Replacing the automorphism $\s$ by the automorphism $f^{-1}\s$, we may assume that $f=e$, and so $\s' = e'\in \mE_{n-1}$. By Lemma \ref{c9Feb12}, $\s'\in \im (\chi_n) \cap \mE_{n-1}= \{ e\}$. Therefore, $I=\sh_{n-2}$.

4. Notice that $\Sh_{n-2}\times \mE_n\subseteq \ker (\res_n)$. By (\ref{FFn=K}), $\CF_n = \mF_n \ltimes \ker (\res_n)$. By statement 1, $\CF_n = \mF_n \times \Sh_{n-2}\times \mE_n$. Therefore, $\ker (\res_n) = \Sh_{n-2} \times \mE_n$ and $\CF_n = \mF_n \times \ker (\res_n)$. $\Box $


It follows from the inclusion $\sh_{n-1}=e^{K\der_{n-1}}\subseteq \mF_{n} =1+\der_{n-1}K[[\der_{n-1}]]$ (see (\ref{Fn=idnt})) that
\begin{equation}\label{Fn=shF}
\mF_n = \sh_{n-1}\times \mF_{n}'
\end{equation}
 where $ \mF_{n}'=1+\der_{n-1}^2 K[[\der_{n-1}]] =1+\d_{n-1}^2K[[\d_{n-1}]]$ (see (\ref{Fn=idnt})). So,
\begin{equation}\label{Fns=def}
\mF_n' = \{ f\in 1+\der_{n-1}^2K[[\der_{n-1}]]\; | \; f(p_i\der_i):=\begin{cases}
p_i\der_i & \text{if }i=1, \ldots , n-1,\\
f(p_n)\der_n & \text{if }i=n,
\end{cases} \\
  {\rm where } \; p_i\in P_{i-1}, i=1, \ldots , n \}.
\end{equation}

 Moreover, $\mF_n=\prod_{i\geq 1} e^{K\d_{n-1}^i}\simeq K^\N$ and $\mF_n'=\prod_{i\geq 2} e^{K\d_{n-1}^i}\simeq K^\N$.

 The next theorem describes the group $G_n$ as an exact product of its explicit subgroups.

\begin{theorem}\label{5Feb12}
Let $\mI := (1+t^2K[[t]], \cdot)$ and $\mJ:= (tK[[t]], +)$. Then for all $n\geq 2$,
\begin{enumerate}
\item $G_n = \mT^n \ltimes (\CT_n \times_{ex} (\Sh_{n-2}\times \mF_n\times \mE_n))=\TAut_K(P_n)_n\times_{ex} (\mF_n'\times \mE_n)$,
\item $ G_n\simeq \TAut_K(P_n)_n\times_{ex} (\mI \times \mJ^{n-2})$.
\end{enumerate}
\end{theorem}

{\it Proof}. 1. The first equality follows from Theorem \ref{26Jan12}.(2) and Theorem \ref{11Feb12}.(1). The second equality follows from the first one, the equality $\mF_n = \sh_{n-1}\times \mF_n'$ (see (\ref{Fn=shF})) and Proposition \ref{c27Jan12}.(3).

2. Statement 2 follows from statement 1 and the facts that $\mF_n'\simeq \mI $  and $\mE_n\simeq \mJ^{n-2}$ (Corollary \ref{d4Feb12}). $\Box $


In Section \ref{GnITER}, Theorem \ref{5Feb12} will be strengthen (Theorem \ref{28Feb12}). Roughly speaking, the exact product will be replaced by a semi-direct product.


\section{The group of automorphism of the   Lie algebra $\ggu_n$ is an iterated semi-direct  product}\label{GnITER}

The aim of this section is to show that the group $G_n$ is an iterated semi-direct  product $\mT^n\ltimes (\UAut_K(P_n)_n\rtimes (\mF_n'\times \mE_n))$ (Theorem \ref{28Feb12}), that none of its subgroups $\TAut_K(P_n)_n$, $\mF_n'\times \mE_n$, $\mF_n\times \mE_n$, $\mF_n'$, $\mF_n$, $\mE_n$ and $\mE_{n,i}$ is a normal subgroup (Corollary \ref{c28Feb12}), and to give characterizations of the groups $\mF_n$, $\mF_n'$ and $\mE_n$ of $G_n$ in invariant terms (Proposition \ref{b12Feb12}). The proof of the results are based on the following  two lemmas.

\begin{lemma}\label{a28Feb12}
Let $t_\l \in \mT^n$ where $\l = (\l_1, \ldots , \l_n) \in K^{*n}$; $e_i'(s_i) \in \mE_{n,i}$ where $ i=2, \ldots , n-1$,  $s_i = \sum_{j\geq 1} \l_{ij}\der_{i-1}^j$ and $\l_{ij}\in K$; and $f=1+\sum_{k\geq 1} \mu_k \d_{n-1}^k\in \mF_n$ where $\d_{n-1}= \ad (\der_{n-1})$. Then
\begin{enumerate}
\item $t_\l e_i'(s_i) t^{-1}_\l = e_i'(\l_i\l_n^{-1}t_\l s_it_\l^{-1})\in \mE_{n,i}$ where $t_\l s_it_\l^{-1}=\sum_{j\geq 1}\l_{ij}\l_{i-1}^{-j}\der_{i-1}^j$.
\item $t_\l f t^{-1}_\l = 1+\sum_{k\geq 1}\mu_k\l^{-k}_{n-1}\d_{n-1}^k\in \mF_n$.
\item If, in addition, $f\in \mF_n'$ (i.e.,  $\mu_1=0$) then $t_\l f t^{-1}_\l = 1+\sum_{k\geq 2}\mu_k\l^{-k}_{n-1}\d_{n-1}^k\in \mF_n'$.
\end{enumerate}
\end{lemma}

{\it Proof}. 1. Let $\s=t_\l e_i'(s_i)t_\l^{-1}$ and $\s':= e_i'(\l_i\l_n^{-1}t_\l s_it_\l^{-1})$. Then $\s (x^\alpha \der_j) = x^\alpha \der_j = \s' (x^\alpha \der_j)$ for all $j\neq i$ and $\alpha \in \N^{j-1}$. For all elements $\alpha \in \N^{i-1}$,
\begin{eqnarray*}
\s: &x^\alpha \der_i& \stackrel{t_\l^{-1}}{\mapsto} \l^{-\alpha}\l_i x^\alpha \der_i \stackrel{e_i'(s_i)}{\mapsto}
t_\l^{-1} (x^\alpha \der_i) +\l^{-\alpha}\l_i s_i(x^\alpha ) \der_n\\
 &\stackrel{t_\l}{\mapsto}& x^\alpha \der_i+\l^{-\alpha}\l_i t_\l s_it_\l^{-1}(t_\l (x^\alpha ))\cdot \l_n^{-1}\der_n = x^\alpha \der_i +\l_i\l_n^{-1}t_\l s_i t_\l^{-1} (x^\alpha ) \der_n=\s'(x^\alpha \der_i).
\end{eqnarray*}
Therefore, $\s = \s'$.

2 and 3. Straightforward. $\Box $


Let $G$ be a group and $a,b\in G$. Then $[a,b]:= aba^{-1}b^{-1}$ is the (group) {\em commutator} of  elements $a$ and $b$ of the group $G$.

\begin{lemma}\label{b28Feb12}
Let $\tau = e^{a_s\der_s}\in \UAut_K(P_n)_n$ where $a_s\in P_{s-1}$ and $1\leq s\leq n$; $e_i'(s_i) \in \mE_{n,i}$ where $ 2\leq i \leq n-1$, $s_i = \sum_{j\geq 1} \l_{ij}\der_{i-1}^j$ and $\l_{ij} \in K$; and $f=1+\sum_{i\geq 1} \mu_i \d_{n-1}^i\in \mF_n$ where $\mu_i\in K$ and $\d_{n-1} = \ad (\der_{n-1})$. Then
\begin{enumerate}
\item $[ e_i'(s_i), \tau ] = \begin{cases}
e^{ s_i(a_i)\der_n}& \text{if }s=i,\\
e& \text{if }s\neq i.
\end{cases}$
\item $[\tau ,f ] = \begin{cases}
e& \text{if }1\leq s<n,\\
e^{-f'(a_n)\der_n}& \text{if }s=n,
\end{cases}$
where $f'(a_n) := \sum_{i\geq 1} \mu_i \der_{n-1}^i(a_n)$.

\end{enumerate}
\end{lemma}

{\it Proof}.  Recall that $\tau : x_s\mapsto x_s+a_s$, $x_i\mapsto x_i$ for $i\neq s$, $\der_j\mapsto \begin{cases}
\der_j-\frac{\der a_s}{\der x_j}\der_s& \text{if }j<s,\\
\der_j& \text{if }j\geq s.\\
\end{cases}$
In the arguments below, the decomposition $\ggu_n=\oplus_{i=1}^n P_{i-1}\der_i$ is often used.

1.  Since $e'(-s_i)=e'(s_i)^{-1}$, the equality in statement 1 is equivalent  to the equality $$[e'(s_i)^{-1} , \tau ] =
\begin{cases}
e^{-s_i(a_i)}\der_n & \text{if }s=i,\\
e& \text{if }s\neq i.\\
\end{cases}$$
Notice that  $[e'(s_i)^{-1} , \tau ] =e'(s_i)^{-1} \cdot \tau e'(s_i) \tau^{-1}$. Let $ e_i'=e_i'(s_i)$ and $\s =  \tau e'_i \tau^{-1}$. We consider three cases separately: $s<i$, $s=i$, and $s>i$.

Case 1: $s<i$. We have to show  that $\s = e_i'$. The automorphisms $\tau^{\pm 1} $ and $e_i'$ respect the vector spaces $V_-:= \oplus_{1\leq j <i}P_{j-1}\der_j$ and $V_+:= \oplus_{i<j \leq n}P_{j-1}\der_j$ (i.e., the vector spaces are invariant under the action of the automorphisms). Moreover, the automorphism $e_i'$ acts as the identity map on both of them, hence so does the automorphism $\s$. In particular, $\s |_{V_-\oplus V_+}=  e_i' |_{V_-\oplus V_+}$. For all elements $\alpha \in \N^{i-1}$,
\begin{eqnarray*}
\s : \,  x^\alpha \der_i &\stackrel{\tau ^{-1}}{\mapsto} & \tau^{-1} (x^\alpha )\der_i \stackrel{e_i'}{\mapsto} \tau^{-1}(x^\alpha \der_i) +s_i\tau^{-1}(x^\alpha)\der_n\\
& \stackrel{\tau }{\mapsto}&  x^\alpha \der_i + \tau s_i\tau^{-1}(x^\alpha )\der_n= x^\alpha \der_i+s_i(x^\alpha ) \der_n =e_i'(x^\alpha \der_i)
 \end{eqnarray*}
since $\tau s_i (x^\alpha ) = s_i\tau (x^\alpha )$ as $s<i$. Therefore, $ \s = e_i'$.

Case 2: $s=i$. We have to show that  $c:= [ e_i'^{-1} (s_i) , \tau ] = e^{- s_i (a_i) \der_n}$.

The automorphisms $\tau$, $e_i'$ and $\xi := e^{- s_i(a_i) \der_n}$ respect the vector space $V_+$, hence so do the automorphisms $\s$ and $c$. Moreover, the automorphisms $e_i'$ and $\xi$ act on $V_+$ as the identity map.  In particular, $c|_{V_+}=\xi |_{V_+}$. Let $1\leq j <i$ and $\alpha \in \N^{j-1}$. Then
\begin{eqnarray*}
 c: \,  x^\alpha \der_j &\stackrel{\tau ^{-1}}{\mapsto} & x^\alpha (\der_j + \frac{\der a_i}{\der x_j}\der_i)
  \stackrel{e_i'}{\mapsto} \tau^{-1}(x^\alpha \der_j) +s_i(x^\alpha \frac{\der a_i}{\der x_j})\der_n\\
& =&  \tau^{-1}(x^\alpha \der_j) + x^\alpha \der_j s_i(a_i)\der_n \stackrel{\tau }{\mapsto} x^\alpha \der_j+x^\alpha  \der_j s_i(a_i) \der_n = x^\alpha (\der_j+\der_j s_i(a_i)\der_n)\\
&\stackrel{e_i'^{-1}}{\mapsto} &  x^\alpha (\der_j+\der_j s_i(a_i)\der_n)=\xi (x^\alpha \der_j),
\end{eqnarray*}
since $\tau (x^\alpha \der_j s_i(a_i))=x^\alpha \der_js_i(a_i)$ as $x^\alpha \der_j s_i(a_i)\in P_{i-1}$.
Finally, for all elements $\alpha \in \N^{i-1}$,
\begin{eqnarray*}
  c: \,  x^\alpha \der_i &\stackrel{\tau ^{-1}}{\mapsto} & x^\alpha \der_i
  \stackrel{e_i'}{\mapsto} x^\alpha \der_i +s_i(x^\alpha )\der_n
  \stackrel{\tau }{\mapsto} x^\alpha \der_i+\tau s_i(x^\alpha ) \der_n \\
  &= &x^\alpha \der_i+s_i(x^\alpha )\der_n =e_i'(x^\alpha \der_i)
\stackrel{e_i'^{-1}}{\mapsto}   x^\alpha \der_i =\xi (x^\alpha \der_i).
\end{eqnarray*}
Therefore, $c=\xi$.

Case 3: $s>i$. We have to show that $\s = e_i'$. The automorphisms $\tau^{\pm 1}$ and $e_i'$ respect the subspace $V=\oplus_{j\neq i} P_{j-1}\der_j$ of $\ggu_n$. Moreover, $e_i'|_V= \id_V$. Therefore, $\s |_V = \id_V=e_i'|_V$. For all elements $ \alpha \in \N^{i-1}$,

\begin{eqnarray*}
 \s : \,  x^\alpha \der_i &\stackrel{\tau ^{-1}}{\mapsto} & \tau^{-1}(x^\alpha )\der_i + \tau^{-1}(x^\alpha )\frac{\der a_s}{\der x_i}\der_s = x^\alpha \der_i + x^\alpha \frac{\der a_s}{\der x_i}\der_s \stackrel{e_i'}{\mapsto} \tau^{-1}(x^\alpha \der_i) +s_i(x^\alpha ) \der_n\\
& \stackrel{\tau }{\mapsto}& x^\alpha \der_i+ s_i(x^\alpha) \der_n =e_i' (x^\alpha \der_i)
\end{eqnarray*}
since $\tau \der_i\tau^{-1} = \tau (\der_i) = \der_i-\frac{\der a_s}{\der x_i}\der_s$ and $\tau^{\pm 1}(x^\alpha ) = x^\alpha$ (since $s>i$ and $\alpha \in \N^{i-1}$).  Therefore, $\s = e_i'$.

2.  Case 1: $1\leq s<n$. In this case, $\tau \der_{n-1}\tau^{-1} = \der_{n-1}$, and so $\tau \d_{n-1} = \d_{n-1}\tau$. This implies that the maps $\tau$ and $ f=1+\sum_{i\geq 1} \mu_i\d_{n-1}^i$ commute.

Case 2: $s=n$. Let $c:= [\tau , f]$. In this case, both automorphisms $\tau$ and $f$ respect the vector space $P_{n-1}\der_n$. Moreover, the automorphism $\tau$ acts as the identity map on it, hence so does the automorphism $c$. Clearly, the automorphism $e^{-f'(a_n)\der_n}$ acts as the identity map on $P_{n-1}\der_n$. Therefore, $c|_{P_{n-1}\der_n}= e^{-f'(a_n)\der_n}|_{P_{n-1}\der_n}$.  Consider the action of the automorphism $c$ on the elements $x^\alpha \der_i$ where $\alpha \in \N^{i-1}$ and $1\leq i <n$,
\begin{eqnarray*}
 c:\,  x^\alpha \der_i &\stackrel{f^{-1}}{\mapsto} & x^\alpha \der_i \stackrel{\tau^{-1}}{\mapsto} x^\alpha (\der_i +\frac{\der a_n}{\der x_i}\der_n)
 \stackrel{f}{\mapsto} \tau^{-1}(x^\alpha \der_i) + f'(x^\alpha \frac{\der a_n}{\der x_i})\der_n=\tau^{-1} ( x^\alpha \der_i) +x^\alpha f'( \frac{\der a_n}{\der x_i})\der_n\\
&\stackrel{\tau}{\mapsto}&  x^\alpha (\der_i+f'( \frac{\der a_n}{\der x_i})\der_n)= x^\alpha (\der_i+\frac{ \der f'( a_n)}{\der x_i} \der_n)=e^{-f'(a_n)\der_n}(x^\alpha \der_i).
 \end{eqnarray*}
Then $c=e^{-f'(a_n)\der_n}$. $\Box $


Statements 1 and 2 of Lemma \ref{b28Feb12} can be rewritten as follows (where $ 2\leq i \leq n-1$)
\begin{equation}\label{esic}
e_i'(s_i) e^{a_s\der_s}e_i'(s_i)^{-1}=
\begin{cases}
e^{ s_i(a_s) \der_n} e^{a_s\der_s}& \text{if }s=i,\\
e^{a_s\der_s}& \text{if }s\neq i,
\end{cases}
\end{equation}
\begin{equation}\label{esic1}
fe^{a_s\der_s}f^{-1}=
\begin{cases}
e^{a_s\der_s}& \text{if }1\leq s<n,\\
e^{(a_n+f'(a_n)) \der_n}=e^{f(a_n\der_n) }=e^{f(a_n) \der_n} & \text{if }s=n.
\end{cases}
\end{equation}
Recall that the map $P_{n-1}\der_n\ra P_{n-1}$, $p\der_n\mapsto p$, is a $\ggu_n$-module isomorphism. Under this isomorphism the action of the element $\der_{n-1}$ on the ideal $P_{n-1}\der_n$ of the Lie algebra $\ggu_n$, which is $\d_{n-1}=\ad (\der_{n-1})$, becomes the partial derivative $\der_{n-1}$ on the polynomial algebra $P_{n-1}$. So, the expression $f(a_n)$ in (\ref{esic1}) makes sense, it simply means $(1+f') (a_n)$.

The next theorem is one of the main results  of the paper.
\begin{theorem}\label{28Feb12}
Let $\mI := (1+t^2K[[t]], \cdot)$ and $\mJ:= (tK[[t]], +)$. Then for all $n\geq 2$,
\begin{enumerate}
\item $G_n = \mT^n \ltimes (\UAut_K(P_n)_n \rtimes ( \mF_n'\times \mE_n))$.
\item $ G_n\simeq \mT^n \ltimes (\UAut_K(P_n)_n \rtimes (\mI \times \mJ^{n-2})$.
    \item $\UAut_K(P_n)_n$ is a normal subgroup of the group $G_n$.
    \item $\CU_n=\UAut_K(P_n)_n \rtimes ( \mF_n'\times \mE_n)$.
\end{enumerate}
\end{theorem}

{\it Proof}. 3. Statement 3 follows from  statement 1.

1. By Proposition \ref{c27Jan12}.(4), $ \UAut_K(P_n)_n = \CT_n\times_{ex} \Sh_{n-1}= \CT_n\times_{ex} (\Sh_{n-2}\times \sh_{n-1})$. By (\ref{Fn=shF}), $\mF_n = \sh_{n-1}\times \mF_n'$. Then, by Theorem \ref{5Feb12}.(1),
\begin{eqnarray*}
 G_n &=& \mT^n \ltimes ( \CT_n\times_{ex} (\Sh_{n-2}\times \sh_{n-1} \times \mF_n'\times \mE_n))=\mT^n \ltimes (\UAut_K(P_n)_n \times_{ex} (\mF_n'\times \mE_n))\\
 &=& \mT^n \ltimes (\UAut_K(P_n)_n\rtimes (\mF_n'\times \mE_n)), \;\; {\rm by }\;\; (\ref{esic}) \; {\rm and }\; (\ref{esic1}).
\end{eqnarray*}

2. Statement 2 follows from statement 1 (see Theorem \ref{5Feb12}.(2)).

 4. Statement 4 follows from the obvious inclusion $\UAut_K(P_n)_n\rtimes (\mF_n'\times \mE_n)\subseteq \CU_n$,  statement 1 and Proposition \ref{17Dec11}.(2).
  $\Box $


\begin{corollary}\label{x28Feb12}

\begin{enumerate}
\item The group $\CT_n':=\{ [0, \ldots , 0,a_n]\, | \, a_n\in \gm_{n-1}\} = e^{\gm_{n-1}\der_n}$ is a normal subgroup of the groups $G_n$ and $\UAut_K(P_n)_n$.
\item $\UAut_K(P_n)_n=\UAut_K(P_{n-1})\ltimes \CT_n'\subseteq G_n$ (this is the equality of subgroups of $G_n$).
\item $G_n/ \CT_n'=\mT^n\ltimes (\UAut_K(P_{n-1})\times \mF_n'\times \mE_n)$.
\end{enumerate}
\end{corollary}

{\it Proof}.  2. Statement 2 is obvious.

1. Statement 1 follows from statement 2, the equalities (\ref{esic}) and (\ref{esic1}), and Theorem \ref{28Feb12}.(1).

3. Statement 3 follows from Theorem \ref{28Feb12}.(1) and statement 2. $\Box $


We say that subgroups $E$ and $F$ of a group $G$ {\em commute} if $ef=fe$ for all elements $e\in E$ and $f\in F$. By (\ref{esic}), the subgroups $\CT_n'$ and $\mE_n$ commute.

\begin{corollary}\label{c28Feb12}

\begin{enumerate}
\item The group $\TAut_K(P_n)_n$ is not a normal subgroup of the group  $G_n$.
\item None of the groups $\mF_n'\times \mE_n$,  $\mF_n\times \mE_n$, $\mF_n'$, $\mF_n$; $\mE_n$ and $\mE_{n,i}$ where $2\leq i \leq n-1$ and  $n\geq 3$ is a normal subgroup of the group $G_n$.
\end{enumerate}
\end{corollary}

{\it Proof}. The corollary follows from Theorem \ref{28Feb12}.(1),  Lemma \ref{a28Feb12} and Lemma \ref{b28Feb12}.  $\Box $


{\bf A formula for multiplication of two elements of the group $G_n$}. Since  the group $G_n$ is the iterated semi-direct product of four groups and one of them is $\UAut_K(P_n)_n$ which is also a semi-direct product of $n$ of its subgroups, the multiplication in the group $G_n$ looks messy. Surprisingly, it looks less messy than one might expect when we change the order in the presentation of an element of the group $G_n$ as a product of four automorphisms. By Theorem \ref{28Feb12}.(1), every element $\s$ of $G_n$ can be written as the product
\begin{equation}\label{sttef}
\s = \tau te'f'
\end{equation}
where $\tau = [a_1, \ldots , a_n]= e^{a_n\der_n}\cdots e^{a_1\der_1}\in \UAut_k(P_n)_n$ and $a_i\in P_{i-1}$ for $i=1, \ldots , n-1$ and $a_n\in \gm_{n-1}$; $t=t_{(\l_1, \ldots , \l_n)}\in \mT^n$; $e'=e_2'(s_2)\cdots e_{n-1}'(s_{n-1})\in \mE_n$ where $s_i=\sum_{j\geq 1} \nu_{ij}\der_{i-1}^j\in \der_{i-1}K[[\der_{i-1}]]$ for $i=2, \ldots , n-1$ and $\nu_{ij}\in K$; and $f'=\prod_{i\geq 2} e^{\mu_i\der_{n-1}^i}=1+\sum_{i\geq 2} \mu_i'\der_{n-1}^i\in \mF_n'=1+\der_{n-1}^2K[[\der_{n-1}]]$. To make  notations simpler and computations more transparent we write $[a_i]$ for $e^{a_i\der_i}$ sometimes.
\begin{equation}\label{efai}
e'f'[a_1, \ldots , a_n](e'f')^{-1} = e^{(f'(a_n)+\sum_{i=2}^{n-1}s_i(a_i))\der_n}[a_1, \ldots , a_{n-1}]=
 e^{((-1+f')(a_n)+\sum_{i=2}^{n-1}s_i(a_i))\der_n}[a_1, \ldots , a_n].
\end{equation}
In more detail, the second equality is obvious. Using (\ref{esic}) and (\ref{esic1}) we obtain the first one:
\begin{eqnarray*}
 e'f'[a_1, \ldots , a_n](e'f')^{-1}& = & e'f'[a_n][a_1, \ldots , a_{n-1}](e'f')^{-1} =
 f'[a_n]f'^{-1}\cdot e'[a_1, \ldots , a_{n-1}]e'^{-1} \\
 &=& e^{f'(a_n) \der_n} \cdot e'[a_{n-1}]e'^{-1}\cdots  e'[a_i]e'^{-1}\cdots e'[a_2]e'^{-1}\cdot a_1\\
 &=& e^{f'(a_n) \der_n} \cdot e^{s_{n-1}(a_{n-1})\der_n}[a_{n-1}]\cdots  e^{s_i(a_i)\der_n}[a_i]\cdots e^{s_2(a_2)\der_n}[a_2]\cdot a_1\\
 &=& e^{(f'(a_n)+\sum_{i=2}^{n-1}s_i(a_i))\der_n}[a_1, \ldots , a_{n-1}].
\end{eqnarray*}
Let $\s_1=\tau_1t_1e'_1f'_1$ be another element of the group $G_n$ where $\tau_1=[b_1, \ldots , b_n]$. Then using (\ref{efai}) we obtain the formula for multiplication in the group $G_n$:
\begin{equation}\label{efai1}
\s\s_1= \tau e^{t(((-1+f')(b_n)+\sum_{i=2}^{n-1}s_i(b_i) )\der_n)}\o_{t}(\tau_1 )\cdot tt_1\cdot \o_{t_1^{-1}}(e')e_1'\cdot \o_{t_1^{-1}}(f')f_1'
\end{equation}
where $\o_{t_1^{-1}}(g)=t_1^{-1}gt_1$.


{\bf Characterizations of the subgroups $\mF_n$, $\mF_n'$ and $\mE_n$}.
By (Corollary 3.12, \cite{Lie-Un-GEN}), 
the ideal $I_{\o^{n-2}+1} = Kx_{n-1}\der_n +\sum_{\alpha \in \N^{n-2}}Kx^\alpha \der_n$ is the least ideal of the Lie algebra $\ggu_n$ which is a faithful $\ggu_{n-1}$-module (with respect to the adjoint action). Hence, its predecessor $I_{\o^{n-2}}=\sum_{\alpha \in \N^{n-2}}Kx^\alpha \der_n$ is the largest ideal of the Lie algebra $\ggu_n$ which is {\em not} a faithful $\ggu_{n-1}$-module. By (Corollary 5.4, \cite{Lie-Un-GEN}), 
the ideal $P_{n-1}\der_n$ is the least ideal $I$ of the Lie algebra $\ggu_n$ such that the Lie factor algebra $\ggu_n/I$ is isomorphic to the Lie algebra $\ggu_{n-1}$. The next proposition gives characterizations of the groups $\mF_n$, $\mF_n'$ and $ \mE_n$.

\begin{proposition}\label{b12Feb12}
Let $n\geq 2$. Then
\begin{enumerate}
\item $\mF_n'=\Fix_{G_n}(\ggu_{n-1}+I_{\o^{n-2}+1})$ where $I_{\o^{n-2}+1} = Kx_{n-1}\der_n +\sum_{\alpha \in \N^{n-2}}Kx^\alpha \der_n$.
    \item  $\mF_n=\Fix_{G_n}(\ggu_{n-1}+I_{\o^{n-2}})$ where $I_{\o^{n-2}}=\sum_{\alpha \in \N^{n-2}}Kx^\alpha \der_n$.
\item $\mE_n = \Fix_{G_n} ( \CD_n+P_{n-1}\der_n)$ where $\CD_n = \sum_{i=1}^n K\der_i$.
\end{enumerate}
\end{proposition}

{\it Proof}. 1. Let $R$ be the RHS of the equality in statement 1. The inclusion $\mF_n'\subseteq R$ is obvious.  Since $\der_1, \ldots , \der_n \in \ggu_{n-1}+I_{\o^{n-2}+1}$, we have the inclusion $ R\subseteq \Fix_{G_n} (\der_1, \ldots , \der_n) = \CF_n$. By Theorem \ref{11Feb12}.(1), $$\CF_n=\Sh_{n-2}\times \mF_n \times \mE_n= \Sh_{n-2}\times (\sh_{n-1}\times \mF_n')\times \mE_n = \Sh_{n-1}\times \mF_n'\times \mE_n.$$ Now, $R= \mF_n'\times (R\cap (\Sh_{n-1}\times \mE_n))$. By looking at the action of the elements of the group $\Sh_{n-1}\times \mE_n$ on the elements $x_1\der_2, x_2\der_3, \ldots , x_{n-1}\der_n$ of the Lie algebra $\ggu_{n-1}+I_{\o^{n-2}+1}$, we conclude  that $R\cap (\Sh_{n-1}\times \mE_n)=R\cap \mE_n$. Every  element of the group $\mE_n$ is uniquely determined by its  action on $\ggu_{n-1}$. Therefore, $R\cap \mE_n  = \{ e\}$, i.e.,  $R=\mF_n'$.

2.  Let $R$ be the RHS of the equality in statement 2.
 Since $\ggu_{n-1}+I_{\o^{n-2}}\subseteq \ggu_{n-1}+I_{\o^{n-2}+1}$, $\mF_n'=\Fix_{G_n} (\ggu_{n-1}+I_{\o^{n-2}+1})\subseteq \Fix_{G_n} (\ggu_{n-1}+I_{\o^{n-2}})=R$, by statement 1. Clearly, $\sh_{n-1}\subseteq R$. Therefore, $\mF_n = \sh_{n-1}\times \mF_n'\subseteq R$.
 Since $\der_1, \ldots , \der_n \in \ggu_{n-1}+I_{\o^{n-2}}$, we have the inclusion $ R\subseteq \Fix_{G_n} (\der_1, \ldots , \der_n) = \CF_n$. By Theorem \ref{11Feb12}.(1),  $\CF_n=\Sh_{n-2}\times \mF_n \times \mE_n$. Now,
  $R=\mF_n \times (R\cap ( \Sh_{n-2}\times \mE_n))$. By looking at the action of the elements of the group $\Sh_{n-2}\times \mE_n$ on the elements $x_1\der_2, x_2\der_3, \ldots , x_{n-2}\der_{n-1}$ of the Lie algebra $\ggu_{n-1}+I_{\o^{n-2}}$, we conclude  that $R\cap (\Sh_{n-2}\times \mE_n) =R\cap \mE_n$. Every element of the group $\mE_n$ is uniquely determined by its  action on $\ggu_{n-1}$. Therefore, $R\cap \mE_n = \{ e\}$, i.e.,  $R=\mF_n$.

  3.  Let $R$ be the RHS of the equality in statement 3. Then $\mE_n\subseteq R$. Clearly, $R\subseteq \Fix_{G_n} (\der_1, \ldots , \der_n) = \CF_n = \Sh_{n-2}\times \mF_n \times \mE_n$. Now, $R= \mE_n \times (R\cap (\Sh_{n-2}\times \mF_n))= \mE_n$ since $R\cap (\Sh_{n-2}\times \mF_n)=\{ e\}$, by looking at the action of the group $\Sh_{n-2}\times \mF_n$ on the ideal $P_{n-1}\der_n$. $\Box $




\section{The canonical decomposition for an  automorphism of the   Lie algebra $\ggu_n$}\label{CDEC}

By Theorem \ref{5Feb12}.(1), every automorphism $\s \in G_n=\mT^n \ltimes (\CT_n\times_{ex}(\Sh_{n-2}\times \mF_n \times \mE_n))$ is the unique product $\s = t \tau s fe'$ where $t\in \mT^n$, $\tau \in \CT_n$, $s\in \Sh_{n-2}$, $ f\in \mF_n$ and $e'\in \mE_n$. This product is called the {\em canonical decomposition}  for the automorphism $\s\in G_n$.  It is a trivial observation that every  automorphism of a Lie algebra is uniquely determined by its action on any generating set for the Lie algebra.  Our goal is to find explicit formulas for the automorphisms $t$, $\tau $, $s$,  $ f$ and $e'$ via the elements $\{ \s (s)\, | \, s\in S_n\}$ where $S_n$ is a certain set of generators for the Lie algebra $\ggu_n$ (Theorem \ref{14Feb12}).

\begin{theorem}\label{14Feb12}
Let $n\geq 2$.
\begin{enumerate}
\item The set $S_n: =\{ \der_1, x_1^j\der_2, \ldots , x_{i-1}^j\der_i, \ldots , x_{n-1}^j\der_n\, | \, j\in \N\}$ is a set of  generators for the Lie algebra $\ggu_n$.
\item Let $\s \in G_n$ and  $\s = t\tau sfe'$ be its canonical decomposition. Below, explicit formulas are given for the automorphisms $t$, $\tau$, $s$, $f$ and $e'$ via the elements $\{ \s (s) \, | \, s\in S_n\}$.
\begin{enumerate}
\item $t=t_{(\l_1, \ldots , \l_n)}$ where $\s (\der_i) =\l_i^{-1}\der_i+\cdots$ for $i=1, \ldots , n$ where the three dots mean smaller terms with respect to the ordering (i.e.,  an element of $\oplus_{j>i}P_{j-1}\der_j$);
\item $\tau :P_n\ra P_n$, $x_i\mapsto x_i'$,  where $x_1'=x_1$ and $ x_i':= \phi_{i-1}\phi_{i-2}\cdots \phi_1 (x_1)$ for $i=2, \ldots , n$; $\phi_i := \sum_{k\geq 0} (-1)^k \frac{x_i'^k}{k!}\der_i'^k$ and $\der_i':= t^{-1}\s (\der_i)$ for $i=1, \ldots , n-1$;
\item $f=1+\sum_{i\geq 1} f_i \d_{n-1}^i$ where $\d_{n-1} = \ad (\der_{n-1})$, $f_i\in K$ and $f_i\der_n = \Phi_{n-1}(t\tau )^{-1} \s (\frac{x_{n-1}^i}{i!}\der_n)$ where $\Phi_{n-1}:= \sum_{k\geq 0} (-1)^k \frac{x_{n-1}^k}{k!}\d_{n-1}^k$ and $\d_{n-1}= \ad (\der_{n-1})$.
\item $s(x_i) =x_i+\mu_i$ for $i=1, \ldots , n-2$ where $(t\tau f)^{-1} \s ( x_i\der_{i+1}) = \mu_i \der_{i+1}+x_i \der_{i+1}+\cdots $(the three dots  denote  an element of $\oplus_{j>i}P_{j-1}\der_j$);
\item by (\ref{E2n1}) and (\ref{E2n}), $e'= e_2'\cdots e_{n-1}'$ is the unique product where $e_i'\in \mE_{n,i}$ for $i=2, \ldots , n-1$, and, for all $j=1, \ldots , n$ and $\alpha \in \N^{j-1}$, $$e_i'(x^\alpha \der_j) =\begin{cases}
x^\alpha\der_i+s_i(x^\alpha ) \der_n& \text{if }j=i,\\
x^\alpha\der_j& \text{if }j\neq i,\\
\end{cases}$$ where $s_i=\sum_{j\geq 1} \nu_{ij}\der_{i-1}^j$, $\nu_{ij}\in K$, $\nu_{ij}\der_n = \Phi_{i-1} ((t\tau sf)^{-1}\s -1) (\frac{x_{i-1}^j}{j!}\der_i)$, $\Phi_{i-1}:= \sum_{k\geq 0} (-1)^k \frac{x_{i-1}^k}{k!}\d_{i-1}^k$ and $\d_{i-1}=\ad (\der_{i-1})$.
\end{enumerate}
\end{enumerate}
\end{theorem}

{\it Proof}. 1. We use induction on $n\geq 2$. The initial case $n=2$ is obvious as the set $S_2$ is a $K$-basis for the Lie algebra $\ggu_2$. Suppose that $n>2$ and the result holds for all $n'<n$. By induction, $S_{n-1}$ is a set of generators for the Lie algebra $\ggu_{n-1}$. Notice that $\ggu_n = \ggu_{n-1}\oplus P_{n-1}\der_n$ and, for all elements $\alpha \in \N^{n-2}$ and $j\in \N$, $x^\alpha x_{n-1}^j \der_n = [ x^\alpha \der_{n-1} , (j+1)^{-1} x_{n-1}^{j+1}\der_n]$. To finish the proof of statement 1 notice that $S_n =S_{n-1}\cup \{ x_{n-1}^j\der_n\, | \, j\in \N\}$ and  the set of  elements $\{ x^\alpha x_{n-1}^j \der_n \}$ is a $K$-basis for the vector space $P_{n-1}\der_n$.

2. Statement (a) is obvious (Proposition \ref{17Dec11}). For all elements $i=1, \ldots , n$, $sfe (\der_i ) =\der_i$ (since $sfe\in \CF_n$, Theorem \ref{11Feb12}.(1)). Then $\tau (\der_i) = \tau sfe (\der_i) = t^{-1}\s (\der_i)$ and statement (b) follows from Theorem \ref{B27Jan12}.(2). The automorphisms $s$ and $e'$ act as the identity map on the vector space $V:= K[x_{n-1}]\der_n$. Therefore, $f|_V=fse'|_V= sfe'|_V= (t\tau )^{-1} \s |_V:V\ra V$ and $f = 1+\sum_{i\geq 1} f_i\d_{n-1}^i$ for some scalars $f_i\in K$. By Proposition \ref{b14Feb12}, $f_i \der_n =\Phi_{n-1}(t\tau )^{-1} \s (\frac{x_{n-1}^i}{i!}\der_n)$. This finishes the proof of statement (c).

For all elements $i=1, \ldots , n-2$, $$(t\tau f)^{-1} \s (x_i\der_{i+1}) = se (x_i\der_{i+1}) = s(x_i\der_{i+1}+\cdots ) = \mu_i\der_{i+1}+x_i\der_{i+1}+\cdots  , $$ and statement  (d) follows.

By (\ref{E2n1}) and (\ref{E2n}), $e'= e_2'\cdots  e_{n-1}'$ is the unique product where $e_i'\in \mE_{n,i}$ for $i=2, \ldots , n-1$, and, for all $j=1, \ldots , n$ and $\alpha \in \N^{j-1}$, $$e_i'(x^\alpha \der_j) =\begin{cases}
x^\alpha\der_i+s_i(x^\alpha ) \der_n& \text{if }j=i,\\
x^\alpha\der_j& \text{if }j\neq i,\\
\end{cases}$$ where $s_i=\sum_{j\geq 1} \nu_{ij}\der_{i-1}^j$ and  $\nu_{ij}\in K$. For each $i=0, \ldots , n-1$,  $$(e-1)|_{K[x_{i-1}]\der_i}=((t\tau sf)^{-1}\s  -1)|_{K[x_{i-1}]\der_i}:
  K[x_{i-1}]\der_i\ra K[x_{i-1}]\der_n, \;\; p\der_i\mapsto s_i(p)\der_n$$ where $p\in K[x_{i-1}]$. By Proposition \ref{b14Feb12},
 $\nu_{ij}\der_n = \Phi_{i-1} ((t\tau sf)^{-1}\s -1) (\frac{x_{i-1}^j}{j!}\der_i)$.  $\Box $



\section{The adjoint group of automorphisms  of the Lie algebra $\ggu_n$}\label{ADJGRP}

The aim of  this section is to show that the adjoint group $\CA (\ggu_n)$ of automorphisms  of the Lie algebra $\ggu_n$ is equal to the group $\UAut_K(P_n)_n$ (Theorem \ref{8Apr12}).

Let $\CG$ be a Lie algebra over the field  $K$ and $\LN (\CG )$ be the set of locally nilpotent elements of the Lie algebra $\CG$. Recall that an element $g\in \CG$ is called a {\em locally nilpotent element} if the inner derivation $\ad (g)$ of the Lie algebra $\CG$ is a locally nilpotent derivation. The set $\LN (\CG )$ is an $\Aut_K(\CG )$-invariant set. Each locally nilpotent element $g$ yields the automorphism $e^{\ad (g)}:=\sum_{i\geq 0} \frac{\ad (g)^i}{i!}$ of the Lie algebra $\CG$ which is called an {\em inner automorphism} of the Lie algebra $\CG$. The subgroup of $\Aut_K(\CG )$, $\CA (\CG ) :=\langle e^{\ad (g)}\, | \, a\in \LN (\CG )\rangle$,  is called the {\em adjoint group} (of automorphisms) of the Lie algebra $\CG$. The adjoint group $\CA (\CG )$ is a {\em normal } subgroup of the group $\Aut_K(\CG)$ since $\s e^{\ad (g)}\s^{-1} = e^{\ad (\s (g))}$ for all  automorphisms $\s\in \Aut_K(\CG )$.


The aim of  this section is to prove the next theorem.

\begin{theorem}\label{8Apr12}

\begin{enumerate}
\item $\CA (\CG ) = \UAut_K(P_n)_n$.
\item The map $\UAut_K(P_n)_n \ra \CA (\CG )$, $ e^a\mapsto e^{\ad (a)}$, is the identity map where $a\in \ggu_n'$ (recall that $\UAut_K(P_n)_n\subset G_n$), i.e., for all elements $u\in \ggu_n$, $e^aue^{-a}= e^{\ad (a)}(u)$.
\end{enumerate}
\end{theorem}
The proof of Theorem \ref{8Apr12}, which is given at the end of the section, is based on the following  proposition that is interesting on its own.

The Lie algebra $\Der_K(P_n)$ is a left $P_n$-module, and so $P_n\ggu_n\subseteq \Der_K(P_n)$. The polynomial algebra $P_n$ is a left $\Der_K(P_n)$-module and a left $\ggu_n$-module. The action of an element $\d \in \Der_K(P_n)$ on the polynomial algebra $P_n$ is denoted either by  $\d *p$ or $\d (p)$. Every element $u\in \ggu_n$ is a locally nilpotent derivation of the polynomial algebra $P_n$ (Proposition \ref{a8Dec11}.(4)). Then $e^u\in \Aut_K(P_n)$.

\begin{proposition}\label{14Apr12}
Let $u,v\in \ggu_n$ and $p\in P_n$. Then
\begin{enumerate}
\item $e^u(v*p)=e^{\ad (u)}(v)*e^u(p)$ (where $ e^u\in \Aut_K(P_n)$).
\item $e^{\ad (u)}(v) = e^u v e^{-u}$ (where $ e^{\ad (u)}\in G_n$).
\item $e^{\ad (u)} (pv)=e^u(p)e^{\ad (u)}(v)$.
\end{enumerate}
\end{proposition}

{\it Proof}. 1.
\begin{eqnarray*}
e^u(v*p)&=& (\sum_{i\geq 0}\frac{u^i}{i!})v(p)=\sum_{i\geq 0}\frac{1}{i!}\sum_{j=0}^i{i\choose j}\ad (u)^j(v)u^{i-j}(p)\\
&=& \sum_{i\geq 0}\sum_{j+k=i}\frac{\ad (u)^j}{j!}(v) \frac{u^k}{k!}(p)=e^{\ad (u)}(v)*e^u(p).  \\
\end{eqnarray*}
2. Recall that $e^{-u}\in \UAut_K(P_n)_n$. In statement 1, replacing the polynomial $p$  by the polynomial $e^{-u}(p)$, we have the equality $e^uve^{-u}*p=e^{\ad (u)}(v)*p$, for all polynomials $p\in P_n$. Therefore, $e^uve^{-u}=e^{\ad (u)}(v)$.

3.  For all natural numbers $s\gg 0$, $\ad (u)^s(pv)=0$,  $u^s (p)=0$ and $\ad (u)^s(v)=0$. So, the infinite sums below are finite sums:
\begin{eqnarray*}
e^{\ad (u)}(pv)&=& (\sum_{i\geq 0}\frac{\ad (u)^i}{i!})(pv)=\sum_{i\geq 0}\frac{1}{i!}\sum_{j=0}^i{i\choose j} u^j(p)\ad (u)^{i-j}(v)\\
&=& \sum_{i\geq 0}\sum_{j+k=i}\frac{u^j}{j!}(p) \frac{\ad (u)^k}{k!}(v)=e^u(p)e^{\ad (u)}(v).  \;\; \Box
\end{eqnarray*}


{\bf Proof of Theorem \ref{8Apr12}}. 1. By Proposition \ref{a28Jan12}, the map $\ggu_n\ra \UAut_K(P_n)$, $ u\mapsto e^u$, is a bijection. In particular, $\UAut_K(P_n) = \{ e^u \, | \, u\in \ggu_n\}$. By Proposition \ref{a27Jan12}.(2), the map (where $v\in \ggu_n$)
$$ \exp : \UAut_K(P_n)\ra G_n , \;\; e^u\mapsto (v\mapsto e^uve^{-u}), $$ is a group homomorphism such that $\ker (\exp )=\sh_n$ and $\im (\exp ) \simeq \UAut_K(P_n) / \sh_n = \UAut_K(P_n)_n$. By Proposition \ref{14Apr12}.(2), $e^uve^{-u}= e^{\ad (u)}(v)$ for all elements $u,v\in \ggu_n$. It follows from this fact that, for all elements $u_1, \ldots , u_s\in \ggu_n$,
$$ e^{\ad (u_1)}\cdots e^{\ad (u_s)}(v) = e^{u_1}\cdots e^{u_s}v(e^{u_1}\cdots e^{u_s})^{-1}. $$
By Proposition \ref{a28Jan12}.(2), $e^{u_1}\cdots e^{u_s}=e^u$ for some element $u\in \ggu_n$. Then, $ e^{\ad (u_1)}\cdots e^{\ad (u_s)}(v) = e^uve^{-u}$, i.e., $ e^{\ad (u_1)}\cdots e^{\ad (u_s)}\in \UAut_K(P_n)_n$.

2. Statement 2 follows from Proposition \ref{14Apr12}.(2). $\Box$

$${\bf Acknowledgements}$$

 The work is partly supported by  the Royal Society  and EPSRC.

\small{

Department of Pure Mathematics

University of Sheffield

Hicks Building

Sheffield S3 7RH

UK

email: v.bavula@sheffield.ac.uk}

\end{document}